\documentclass[final,leqno]{siamltex}
\usepackage{amsmath,amsfonts,amssymb}

\usepackage{graphicx}
\usepackage{mathtools}

\usepackage{tikz}
\usetikzlibrary{decorations.pathreplacing}

\newtheorem{assumption}{Assumption}[section]

\def\0{\phantom{0}}
\def\d{\partial}
\def\twoVec(#1,#2){\left(\begin{matrix}#1\\#2\end{matrix}\right)}
\def\grad{\nabla}
\def\div{\nabla\cdot}
\def\Div{\textrm{div}}
\def\Curl{\textrm{curl}}

\def\spn{\textrm{span}}

\def\Bu{\mathbb B}
\def\Po{\mathbb P}
\def\Supp{\mathbb S}
\def\Re{\mathbb R}

\def\F{{\mathbf F}}
\def\V{{\mathbf V}}

\def\v{{\mathbf v}}
\def\u{{\mathbf u}}
\def\x{{\mathbf x}}
\def\y{{\mathbf y}}
\def\bfpsi{\pmb\psi}

\def\cA{{\cal A}}

\def\cDS{{\cal{DS}}}

\def\cL{{\cal L}}
\def\cN{{\cal N}}

\def\cP{{\cal P}}
\def\cS{{\cal S}}
\def\cT{{\cal T}}

\def\Line(#1,#2)(#3,#4){\qbezier(#1,#2)(#1,#2)(#3,#4)}
\def\longlongrightarrow%
{\DOTSB\protect\relbar\protect\joinrel\protect\relbar\protect\joinrel\protect\relbar\protect\joinrel\rightarrow}
\def\hooklongrightarrow{\DOTSB\lhook\joinrel\longrightarrow}

\title{Direct Serendipity and Mixed Finite Elements\\on Convex Polygons%
\thanks{This work was supported by the U.S.~National Science Foundation under grant DMS-2111159.}}

\author{
  Todd Arbogast\thanks{University of Texas at Austin;
  Department of Mathematics, C1200; Austin, TX 78712-1202 and
  Oden Institute for Computational Engineering and Sciences, C0200;
  Austin, TX 78712-1229 (ORCID ID 0000-0001-9692-5478, arbogast@oden.utexas.edu)}
\and
 Chuning Wang\thanks{University of Texas at Austin;
 Department of Mathematics, C1200; Austin, TX 78712-1202 
 (cwangaw@utexas.edu)}
}

\begin{document}

\maketitle

\date{\today}

\begin{abstract}
  We construct new families of \emph{direct} serendipity and \emph{direct} mixed finite elements on
  general planer convex polygons that are $H^1$ and $H(\Div)$ conforming, respectively, and possess
  optimal order of accuracy for any order.  They have a minimal number of degrees of freedom subject
  to the conformity and accuracy constraints. The name arises because the shape functions are
  defined \emph{directly} on the physical elements, i.e., without using a mapping from a reference
  element. The finite element shape functions are defined to be the full spaces of scalar or vector
  polynomials plus a space of supplemental functions.  The direct serendipity elements are the
  precursors of the direct mixed elements in a de Rham complex. The convergence properties of the
  finite elements are shown under a regularity assumption on the shapes of the polygons in the mesh,
  as well as some mild restrictions on the choices one can make in the construction of the
  supplemental functions. Numerical experiments on various meshes exhibit the performance of these
  new families of finite elements.
\end{abstract}




\pagestyle{myheadings} \thispagestyle{plain}
\markboth{Todd Arbogast and Chuning Wang}
{Direct Serendipity and Mixed Finite Elements}


\section{Introduction}\label{sec:Intro}

Serendipity finite elements defined on a rectangle $\hat{E}$, denoted
as $\cS_{r}(\hat{E})$, $r\ge1$, are well known to be $H^1$-conforming
and approximate to order $r+1$ with a minimal number of degrees of
freedom (DoFs). The finite elements $\cS_{r+1}(\hat E)$ are related to Brezzi-Douglas-Marini \cite{BDM_1985} mixed finite elements $\textrm{BDM}_{r}(\hat{E})$, $r\ge1$, through a de Rham complex \cite{AFW_2010_feec}. $\textrm{BDM}_{r}(\hat{E})$ is $H(div)$-conforming and has optimal order approximation properties with a minimal number of DoFs. Arnold and Awanou \cite{Arnold_Awanou_2011,Arnold_Awanou_2014} have given a definition, construction, and geometric decomposition of $\cS_{r}(\hat{E})$ and $\textrm{BDM}_{r}(\hat{E})$ of any approximation order on cubical meshes in any dimension. However, the elements lose optimal order accuracy when mapped to a quadrilateral~$E$.

Recently, the current authors and Z.\ Tao \cite{Arbogast_Tao_Wang_2020x_serendipityMixed} constructed serendipity spaces directly on quadrilaterals of any approximation order $r+1\ge2$ without using a mapping from a reference element. The resulting new family of spaces were called \emph{direct} serendipity finite elements and denoted $\cDS_r(E)$, $r\ge1$. The de Rham complex then yields a strategy to construct $H(\Div)$ conforming \emph{direct} mixed finite elements, denoted $\V_r^{r-1}(E)$ and $\V_r^r(E)$, giving optimal order reduced and full $H(\Div)$-approximation with a minimal number of DoFs. The direct serendipity finite elements take the form
\begin{equation}\label{eq:cDS=P+S_quad}
	\cDS_r(E) = \Po_r(E)\oplus\Supp_r^\cDS(E),
\end{equation}
where $\Po_r(E)$ is the space of polynomials on $E$ up to degree $r$, and $\Supp_r^\cDS(E)$ consists of supplemental functions. The direct mixed elements take a similar form.
In this paper, we construct a new family of direct serendipity and direct mixed finite elements for a general planar convex polygon, discuss their approximation properties, and test their performance by numerical experiments.

Other approaches to construct serendipity and mixed finite elements with a minimal number of degrees of freedom have appeared in the literature. In \cite{Rand_Gillette_Bajaj_2014}, Rand, Gillette, and Bajaj used products of linear generalized barycentric coordinates to construct serendipity finite elements on quadrilaterals. Based on this work, Sukumar \cite{Sukumar_2013} constructed quadratic maximum-entropy serendipity shape functions. These two works only have elements with quadratic order of accuracy, and it appears to be technically difficult to develop higher order accurate serendipity finite elements in this way. However, their construction works for general polygonal elements, including non-convex ones. For mixed spaces, Chen and Wang \cite{Chen_Wang_2017} constructed minimal degree $H(\Curl)$ and $H(\Div)$ conforming finite elements of linear accuracy based on generalized barycentric coordinates and the Whitney forms. Floater and Lai \cite{Floater_Lai_2016} generalized this idea to construct finite element spaces for a general order of accuracy~$r$. However, their construction asks for more DoFs than the minimum, since $\frac{1}{2}(r-1)(r-2)$ interior DoFs are always required for any polygon. Another methodology, the serendipity virtual element method, was introduced in \cite{BBMR_2016_serendipityVEM} to deal with general polygonal elements, including non-convex and very distorted elements. The method works for any order of accuracy $r$, but it uses even more interior DoFs, $\frac{1}{2}r(r-1)$.

In the rest of this paper, we generalize the construction in \cite{Arbogast_Tao_Wang_2020x_serendipityMixed} to a general convex polygon $E_N$ with $N$ vertices. We begin by introducing some notation in Section~\ref{sec:notation}. In Section~\ref{sec:DShigh} we define higher order direct serendipity elements ($r \ge N-2$) and show their unisolvence and conformity by the construction of nodal basis functions. In Section~\ref{sec:DSlow}, lower order direct serendipity elements ($r < N-2$) are constructed within a higher order direct serendipity space. We discuss the approximation properties and convergence rates of the space $\cDS_r$ over the whole domain $\Omega$ in Section~\ref{sec:approxDS}. In Sections~\ref{sec:deRham} and \ref{sec:approxMixed}, we construct direct mixed finite elements from the direct serendipity elements and the de Rham complex, and then discuss the convergence theory. In Section~\ref{sec:numerics}, we provide some numerical results that test the performance of our direct spaces on various meshes. Finally, the results are summarized in Section~\ref{sec:conc}.


\section{Some notation}\label{sec:notation}

Let $\Po_r(\omega)$ denote the space of polynomials of degree up to $r$ on $\omega\subset\Re^d$,
where $d=0$ (a point), $1$, or~$2$. Recall that
\begin{equation}\label{eq:dimPo}
\dim\Po_r(\Re^d) = \twoVec(r+d,d) = \frac{(r+d)!}{r!\,d!}.
\end{equation}
Let $\tilde\Po_r(\omega)$ denote the space of homogeneous polynomials of degree $r$ on $\omega$.  Then
\begin{equation}\label{eq:dimTildePo}
\dim\tilde\Po_r(\Re^d) = \twoVec(r+d-1,d-1) = \frac{(r+d-1)!}{r!\,(d-1)!},\quad d\ge1.
\end{equation}

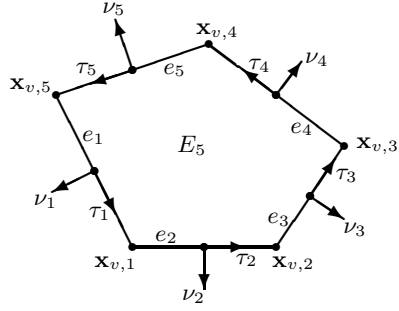
\begin{figure}[ht]\centering
\setlength\unitlength{3.2pt}
\begin{picture}(46,32)(-6,-6)\small
%
\thicklines
\put(3,-2){\line(1,0){17}}
\put(20,-2){\line(2,3){8}} 
\put(28,10){\line(-4,3){16}}
\put(12,22){\line(-3,-1){18}}
\put(-6,16){\line(1,-2){9}}

\put(11.5,-2){\circle*{1}}\put(11.5,-2){\vector(0,-1){5}}
\put(24,4){\circle*{1}}\put(24,4){\vector(3,-2){4}}
\put(20,16){\circle*{1}}\put(20,16){\vector(3,4){3}}
\put(3,19){\circle*{1}}\put(3,19){\vector(-1,3){2}}
\put(-1.5,7){\circle*{1}}\put(-1.5,7){\vector(-2,-1){5}}

\put(10.2,-7.8){\makebox(0,0){$\nu_2$}}
\put(29.3,0.1){\makebox(0,0){$\nu_3$}}
\put(24.8,20.2){\makebox(0,0){$\nu_4$}}
\put(0.8,26.2){\makebox(0,0){$\nu_5$}}
\put(-7.3,3.4){\makebox(0,0){$\nu_1$}}

\put(11.5,-2){\circle*{1}}\put(11.5,-2){\vector(1,0){5}}
\put(24,4){\circle*{1}}\put(24,4){\vector(2,3){3}}
\put(20,16){\circle*{1}}\put(20,16){\vector(-4,3){4}}
\put(3,19){\circle*{1}}\put(3,19){\vector(-3,-1){5}}
\put(-1.5,7){\circle*{1}}\put(-1.5,7){\vector(1,-2){2.5}}

\put(16.5,-3.5){\makebox(0,0){$\tau_2$}}
\put(28.3,6.5){\makebox(0,0){$\tau_3$}}
\put(18,19.5){\makebox(0,0){$\tau_4$}}
\put(-2.4,19.1){\makebox(0,0){$\tau_5$}}
\put(-0.9,2){\makebox(0,0){$\tau_1$}}

\put(20,-2){\circle*{1}}\put(22,-4){\makebox(0,0){$\x_{v,2}$}}
\put(28,10){\circle*{1}}\put(32,10){\makebox(0,0){$\x_{v,3}$}}
\put(12,22){\circle*{1}}\put(13,23.5){\makebox(0,0){$\x_{v,4}$}}
\put(-6,16){\circle*{1}}\put(-9,17){\makebox(0,0){$\x_{v,5}$}}
\put(3,-2){\circle*{1}}\put(1,-4){\makebox(0,0){$\x_{v,1}$}}

\put(10,10){\makebox(0,0){$E_5$}}

\put(7,-0.8){\makebox(0,0){$e_2$}}
\put(20.3,1.3){\makebox(0,0){$e_3$}}
\put(23,12){\makebox(0,0){$e_4$}}
\put(8,19){\makebox(0,0){$e_5$}}
\put(-1.6,11){\makebox(0,0){$e_1$}}

\end{picture}
\caption{A pentagon $E_5$, with edges $e_i$, outer unit normals $\nu_i$, tangents $\tau_i$, and vertices $\x_{v,i}$.
\label{fig:numbering}}
\end{figure}

Let the element $E=E_N\subset\Re^2$ be a closed, nondegenerate, convex polygon with $N\ge3$ edges.  By nondegenerate, we mean that $E_N$ does not degenerate to any polygon with fewer edges, a line segment, or a point.  We choose to identify the edges and vertices of $E_N$ adjacently in the counterclockwise direction, as depicted in Figure~\ref{fig:numbering} (throughout the paper, we interpret indices modulo~$N$). Let the edges of $E_N$ be denoted $e_i$, $i=1,2,\ldots,N$, and the vertices be $\x_{v,i}=e_i\cap e_{i+1}$. Let $\nu_i$ denote the unit outer normal to edge $e_i$, and let $\tau_i$ denote the unit tangent vector of $e_i$ oriented in the counterclockwise direction, for $i=1,2,\ldots,N$.

Let the overall domain $\Omega\subset\Re^2$ be a connected, polygonal open set with a Lipschitz boundary (i.e., $\Omega$ has no slits).  Let $\cT_h$ be a conforming finite element partition or mesh of $\bar\Omega$ into elements (closed, nondegenerate, convex polygons) of maximal diameter $h>0$.  These elements need not have the same number of edges.

For any two distinct points $\y_1$ and $\y_2$, let $\cL{[\y_1,\y_2]}$ be the line passing through $\y_1$ and $\y_2$, and take $\nu{[\y_1,\y_2]}$ to be the unit vector normal to this line interpreted as going from $\y_1$ to $\y_2$ in the clockwise direction (i.e., pointing to the right). Then we define a linear polynomial giving the signed distance of $\x$ to $\cL{[\y_1,\y_2]}$ as  
\begin{equation}\label{eq:any_two_point_lambda}
\lambda{[\y_1,\y_2]}(\x)= -(\x-\y_2)\cdot\nu{[\y_1,\y_2]}.
\end{equation}
To simplify the notation for linear functions that will be used throughout the paper, let $\cL_i=\cL[\x_{v,i-1},\x_{v,i}]$ be the line containing edge $e_i$ and let $\lambda_i(\x)$ give the distance of $\x\in\Re^2$ to edge~$e_i$ opposite the normal direction, i.e.,
\begin{align}
\label{eq:lambda_i}
\lambda_i(\x) &= \lambda[\x_{v,i-1},\x_{v,i}](\x)= - (\x-\x_{v,i})\cdot\nu_i, \quad i=1,2,\ldots,N.
\end{align}
These functions are strictly positive in the interior of $E_N$, and each vanishes on the edge which defines it.

Recall Ciarlet's definition~\cite{Ciarlet_1978} of a finite element.
\begin{definition}[Ciarlet 1978]\label{defn:ciarlet}
	Let 
	\begin{enumerate}
		\item[$1.$] $E\subset \Re^d$ be a bounded closed set with nonempty interior and a Lipschitz continuous
		boundary,
		\item[$2.$] $\cP$ be a finite-dimensional space of functions on $E$, and
		\item[$3.$] $\cN = \{ N_1, N_2,\ldots, N_{\dim\cP} \}$ be a basis for $\cP'$.
	\end{enumerate}
	Then $(E, \cP, \cN)$ is called a \emph{finite element.}
\end{definition}

\section{Direct serendipity elements when $r\ge N-2$}\label{sec:DShigh}

We construct direct serendipity elements for $r \ge N-2$ in this section.  The construction for $1\le r < N-2$ is different, and it is discussed in Section~\ref{sec:DSlow}.

\begin{table}[ht]
\caption{
Geometric decomposition and number of degrees of freedom (DoFs) associated to each geometric object of a polygon $E_N$ for a serendipity element of index $r\ge N-2\ge1$.\label{tab:geometricDecomp}}
\centerline{\begin{tabular}{ccccc}
Dimension & Object  & Object  & DoFs per  & Total  \\
                 &  Name &  Count &   Object &  DoFs \\
\hline
0 & vertex & $N$ & 1 & $N$\\
1 & edge & $N$ & $\dim\Po_{r-2}(\Re)$ & $N(r-1)$ \\
2 & cell & 1 & $\dim\Po_{r-N}(\Re^2)$ & $\frac12(r-N+2)(r-N+1)$\\
\end{tabular}}
\end{table}

To obtain both that $\Po_r(E)\subset\cDS_r(E)$ and that the shape functions on adjoining elements
can be merged together continuously, we consider the lower dimensional geometric objects within $E$.
As shown in Table~\ref{tab:geometricDecomp}, the minimal number of DoFs associated to each lower
dimensional object must correspond to the dimension of the polynomials that restrict to that object.
A polygon with $N$ sides has $N$ vertices, $N$ edges, and one cell of dimension 0, 1, and 2,
respectively. Each vertex requires $\dim\Po_r(\Re^0)=1$ DoF, each edge requires
$\dim\Po_{r-2}(\Re)=r-1$ DoFs (interior to the edge), and each cell requires
$\dim\Po_{r-N}(\Re^2)=\twoVec(r-N+2,2)=\frac12(r-N+2)(r-N+1)$ DoFs (interior to the cell).  There
are cell DoFs only if $r\ge N$, but the formula works for $r\ge N-2$.  The total number of DoFs is
then $D_{N,r}$, where \begin{equation}\label{eq:amtotalNumber} D_{N,r} = N + N(r-1) +
  \frac12(r-N+2)(r-N+1) = \dim\Po_r(E) + \frac12N(N-3),
\end{equation}
and so to define $\cDS_r(E)$, we will supplement $\Po_r(E)\subset\cDS_r(E)$ with the span of
$\frac12N(N-3)$ linearly independent functions.  The quantity $\frac12N(N-3)$ can be interpreted as
the number of pairs of edges that are not adjacent.


\subsection{Shape functions}\label{sec:supplements}

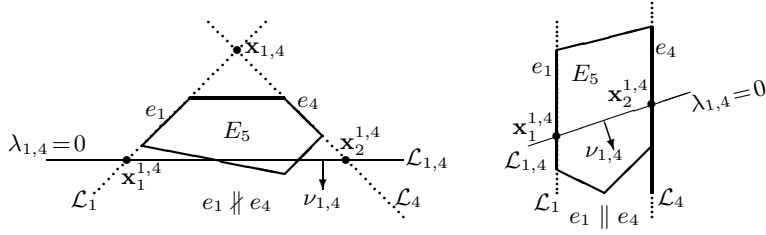
\begin{figure}
	\centering{\setlength\unitlength{1.8pt}\begin{picture}(70,45)(9,2)\small
		%
		\thicklines
		\multiput(10,10)(1.0,1.0){37}{\circle*{0.2}}
		\multiput(75,5)(-1.0,1.0){42}{\circle*{0.2}}
		\put(20,20){\line(1,1){10}}
		\put(58,22){\line(-1,1){8}}
		\put(20,20){\line(5,-1){30}}
		\put(50,14){\line(1,1){8}}
		\put(30,30){\line(1,0){20}}
		\put(40,40){\circle*{1.5}}\put(46,40){\makebox(0,0){$\x_{1,4}$}}
		\put(23,27){\makebox(0,0){$e_1$}}
		\put(80,17){\makebox(0,0){$\cL_{1,4}$}}
		\put(8,8.5){\makebox(0,0){$\cL_1$}}
		\put(55,29){\makebox(0,0){$e_4$}}
		\put(76,9){\makebox(0,0){$\cL_4$}}
		\put(40,23){\makebox(0,0){$E_5$}}
		\thinlines
		\put(0,17){\line(1,0){75}}
		\put(0,20){\makebox(0,0){$\lambda_{1,4}\!=\!0$}}
		\put(17.7,17){\makebox(0,0){\circle*{1.5}}}\put(20,13.5){\makebox(0,0){$\x^{1,4}_1$}}
		\put(63.7,17){\makebox(0,0){\circle*{1.5}}}\put(66,21){\makebox(0,0){$\x^{1,4}_2$}}
		\put(58,17){\vector(0,-1){6.3}}\put(58,8){\makebox(0,0){$\nu_{1,4}$}}
		\put(40,8){\makebox(0,0){$e_1\nparallel e_4$}}
		\end{picture}\qquad\qquad\begin{picture}(50,50)(4,2)\small
		%
		\thicklines
		\multiput(10,10)(0,1.0){37}{\circle*{0.2}}
		\multiput(30,5)(0,1.0){46}{\circle*{0.2}}
		\put(10,15){\line(0,1){25}}
		\put(30,10){\line(0,1){35}}
		\put(10,15){\line(2,-1){10}}
		\put(20,10){\line(1,1){10}}
		\put(10,40){\line(4,1){20}}
		\put(7,36){\makebox(0,0){$e_1$}}
		\put(3.5,16.5){\makebox(0,0){$\cL_{1,4}$}}
		\put(8,8){\makebox(0,0){$\cL_1$}}
		\put(33,40){\makebox(0,0){$e_4$}}
		\put(34,9){\makebox(0,0){$\cL_4$}}
		\put(16,35){\makebox(0,0){$E_5$}}
		\thinlines
		\put(4,20){\line(3,1){34}}
		\put(46,30){\rotatebox{13}{\makebox(0,0){$\lambda_{1,4}\!=\!0$}}}
		\put(10.75,22){\makebox(0,0){\circle*{1.5}}}\put(5,24){\makebox(0,0){$\x^{1,4}_1$}}
		\put(30.75,28.7){\makebox(0,0){\circle*{1.5}}}\put(25,31){\makebox(0,0){$\x^{1,4}_2$}}
		\put(20,25){\vector(1,-3){2}}\put(20,18){\makebox(0,0){$\nu_{1,4}$}}
		\put(20,5){\makebox(0,0){$e_1\parallel e_4$}}
		\end{picture}}
	\caption{Illustration on $E_5$ of the zero line $\cL_{1,4}$ of $\lambda_{1,4}(\x)=-(\x-\x_2^{1,4})\cdot\nu_{1,4}$
		and the intersection point $\x_{1,4}=\cL_1\cap\cL_4$, if it exists.\label{fig:lambda-h}}
\end{figure}

To define the supplemental basis functions, we have two series of choices for each $i,j$ such that
$1\le i<j\le N$ and $2\le j-i\le N-2$ (i.e., $i$ and $j$ are nonadjacent).  First, as shown in
Fig.~\ref{fig:lambda-h}, one must choose two distinct points $\x^{i,j}_1\in\cL_i$ and
$\x^{i,j}_2\in\cL_j$ that avoid the intersection point $\x_{i,j}=\cL_i\cap\cL_j$, if it exists. Then let
\begin{equation}
  \lambda_{i,j}(\x) = \lambda[\x^{i,j}_1, \x^{i,j}_2](\x) = -(\x - \x^{i,j}_2)\cdot\nu_{i,j},
  \quad\nu_{i,j}=\nu[\x^{i,j}_1, \x^{i,j}_2],
\end{equation}
be the linear function associated to the line $\cL_{i,j}=\cL[\x^{i,j}_1,\x^{i,j}_2]$.  Simple
choices are to take the midpoints of the edges, or
\begin{equation}\label{eq:simpleLambdaHV}
  \lambda_{i,j}^{\text{simple}} = \frac{\lambda{[\x_{v,j},\x_{v,i-1}]} - \lambda{[\x_{v,i},\x_{v,j-1}]}}
  {\|\nu{[\x_{v,j},\x_{v,i-1}]} - \nu{[\x_{v,i},\x_{v,j-1}]}\|},
\end{equation}
although the normalization is not strictly necessary.

Second, one must choose the functions $R_{i,j}$ to satisfy the properties
\begin{equation}\label{eq:R-2sided}
\begin{alignedat}3
R_{i,j}(\x)|_{e_i} &= -1,&&\quad&R_{i,j}(\x)|_{e_j} &= 1.
\end{alignedat}
\end{equation}
These are $\pm1$ on $e_i$ and $e_j$, but arbitrary on the other edges. For example, take the simple
rational functions
\begin{equation}\label{eq:R-simple}
R_{i,j}(\x) = R_{i,j}^{\text{simple}}(\x) = \frac{\lambda_i(\x) - \lambda_j(\x)}{\lambda_i(\x) + \lambda_j(\x)}
\end{equation}
(note that the denominators do not vanish on $E_N$, since $e_i$ and $e_j$ are not adjacent). 

The supplemental basis functions are then constructed as
\begin{equation}\label{eq:supplementFcns}
\phi_{s,i,j} = \Big(\prod_{k\neq i,j}\lambda_k\Big)\lambda_{i,j}^{r-N+2}R_{i,j},
\end{equation}
and the supplemental space is defined to be
\begin{align}
\label{eq:supplementSpace}
\Supp_r^{\cDS}(E_N) &= \Supp_r^{\cDS}(E_N;\lambda_{i,j},R_{i,j})\\\nonumber
  &= \spn\big\{\phi_{s,i,j}  : 1\le i<j\le N,\ 2\le j-i\le N-2\big\}.
\end{align}
The $\lambda_{i,j}$'s are not needed when $r=N-2$, and $\Supp_r^{\cDS}(E_N)$ is empty when $N=3$. The
full space $\cP$ in Definition~\ref{defn:ciarlet} is
\begin{equation}\label{eq:cDS=P+S}
\cDS_r(E_N) = \Po_r(E_N)\oplus\Supp_r^\cDS(E_N).
\end{equation}
Each of our earlier choices gives rise to a distinct family of direct serendipity elements of index $r\geq N-2\ge1$.


\subsection{Degrees of Freedom}\label{sec:dofs}

DoFs could be defined in various ways. DoFs based on orthogonal polynomials are generally more
numerically stable.  However, to ease the exposition and proof of unisolvence, we simply use DoF
functionals given by evaluation at (nodal) points.

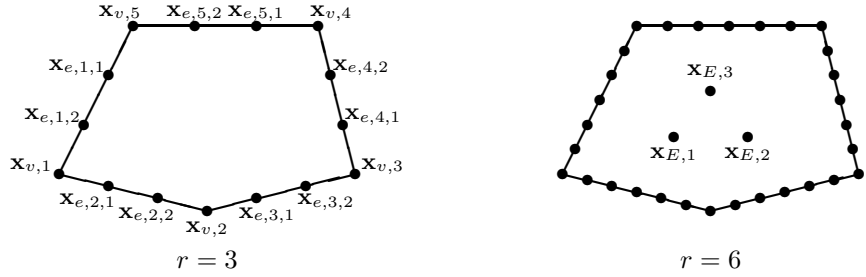
\begin{figure}[htbp]{\setlength\unitlength{3.5pt}
		\centerline{\begin{picture}(43,29.5)(-13.5,-7.5)\small
				\thicklines
				\put(8,-6.25){\makebox(0,0){\normalsize$r=3$}}
				\put(-8,3){\line(1,2){8}}
				\put(-8,3){\line(4,-1){16}}
				\put(8,-1){\line(4,1){16}}
				\put(24,3){\line(-1,4){4}}
				\put(0,19){\line(1,0){20}}
				\put(-8,3){\circle*{1}}
				\put(8,-1){\circle*{1}}
				\put(24,3){\circle*{1}}
				\put(20,19){\circle*{1}}
				\put(0,19){\circle*{1}}
				\put(-2.67,13.67){\circle*{1}}
				\put(-5.33,8.33){\circle*{1}}
				\put(-2.67,1.67){\circle*{1}}
				\put(2.67,0.33){\circle*{1}}
				\put(13.33,0.33){\circle*{1}}
				\put(18.67,1.67){\circle*{1}}
				\put(22.67,8.33){\circle*{1}}
				\put(21.33,13.67){\circle*{1}}
				\put(13.33,19){\circle*{1}}
				\put(6.67,19){\circle*{1}}
				\put(-11,3.8){\makebox(0,0){$\x_{v,1}$}}
				\put(8,-3){\makebox(0,0){$\x_{v,2}$}}
				\put(27,3.8){\makebox(0,0){$\x_{v,3}$}}
				\put(21.4,20.4){\makebox(0,0){$\x_{v,4}$}}
				\put(-1.4,20.4){\makebox(0,0){$\x_{v,5}$}}
				\put(-6,14.5){\makebox(0,0){$\x_{e,1,1}$}}
				\put(-8.67,9.17){\makebox(0,0){$\x_{e,1,2}$}}
				\put(-4.87,0){\makebox(0,0){$\x_{e,2,1}$}}
				\put(1.5,-1.5){\makebox(0,0){$\x_{e,2,2}$}}
				\put(14.5,-1.5){\makebox(0,0){$\x_{e,3,1}$}}
				\put(20.87,0){\makebox(0,0){$\x_{e,3,2}$}}
				\put(25.97,9.17){\makebox(0,0){$\x_{e,4,1}$}}
				\put(24.63,14.5){\makebox(0,0){$\x_{e,4,2}$}}
				\put(13.33,20.4){\makebox(0,0){$\x_{e,5,1}$}}
				\put(6.67,20.4){\makebox(0,0){$\x_{e,5,2}$}}
			\end{picture}\qquad \qquad
			\begin{picture}(43,29.5)(-13.5,-7.5)\small
				\thicklines
				\put(8,-6.25){\makebox(0,0){\normalsize$r=6$}}
				\put(-8,3){\line(1,2){8}}
				\put(-8,3){\line(4,-1){16}}
				\put(8,-1){\line(4,1){16}}
				\put(24,3){\line(-1,4){4}}
				\put(0,19){\line(1,0){20}}
				\put(-8,3){\circle*{1}}
				\put(8,-1){\circle*{1}}
				\put(24,3){\circle*{1}}
				\put(20,19){\circle*{1}}
				\put(0,19){\circle*{1}}
				\put(-1.335,16.335){\circle*{1}}
				\put(-2.67,13.67){\circle*{1}}
				\put(-4,11){\circle*{1}}
				\put(-5.33,8.33){\circle*{1}}
				\put(-6.665,5.665){\circle*{1}}
				\put(-5.335,2.335){\circle*{1}}
				\put(-2.67,1.67){\circle*{1}}
				\put(0,1){\circle*{1}}
				\put(2.67,0.33){\circle*{1}}
				\put(5.335,-0.33){\circle*{1}}
				\put(10.665,-0.33){\circle*{1}}
				\put(13.33,0.33){\circle*{1}}
				\put(16,1){\circle*{1}}
				\put(18.67,1.67){\circle*{1}}
				\put(21.335,2.33){\circle*{1}}
				\put(23.335,5.665){\circle*{1}}
				\put(22.67,8.33){\circle*{1}}
				\put(22,11){\circle*{1}}
				\put(21.33,13.67){\circle*{1}}
				\put(20.665,16.335){\circle*{1}}
				\put(16.665,19){\circle*{1}}
				\put(13.33,19){\circle*{1}}
				\put(10,19){\circle*{1}}
				\put(6.67,19){\circle*{1}}
				\put(3.335,19){\circle*{1}}
				\put(8,12){\circle*{1}}
				\put(4,7){\circle*{1}}
				\put(12,7){\circle*{1}}
				\put(8,14){\makebox(0,0){$\x_{E,3}$}}
				\put(4,5){\makebox(0,0){$\x_{E,1}$}}
				\put(12,5){\makebox(0,0){$\x_{E,2}$}}
		\end{picture}}
		\caption{The nodal points for the DoFs of a direct serendipity finite element $E_5$, for small~$r$.
			\label{fig:nodalDoFs}}
	}
\end{figure}

As depicted in Figure~\ref{fig:nodalDoFs}, for vertex DoFs, the nodal points are exactly the
vertices $\x_{v,i}$, of~$E_N$, where $i=1,2,\ldots, N$.  For edge DoFs, we simply fix nodal
points so that they, plus the two vertices, are equally distributed on each edge. There are $r-1$
nodal points on the interior of each edge, which can be denoted $\x_{e,i,j}$, $j=1,2,\ldots,r-1$,
for nodal points that lie on edge $e_i$, $i=1,2,\ldots, N$, ordered in the counterclockwise
direction. The interior cell DoFs can be set, for example, on points of a triangle~$T$ strictly
inside $E$, where the set of nodal points is the same as the nodes of the Lagrange element of order
$r-N$ on the triangle~$T$.  We denote the interior nodal points as $\x_{E,i}$,
$i=1,2,\ldots,\frac12(r-N+2)(r-N+1)$.

The total number of nodal points is indeed $D_{N,r}$.  If
$\{x_1^{\textrm{nodal}},x_2^{\textrm{nodal}},\ldots,x_{D_{N,r}}^{\textrm{nodal}}\}$ is the set of
all nodal points, then the set of DOFs ($\cN$ in Definition~\ref{defn:ciarlet}) is
\begin{equation}\label{eq:nodalFunctionals}
\cN = \{N_i : N_i(\phi) = \phi(\x^{\textrm{node}}_i)\text{ for all }\phi(\x),\ i=1,2,\ldots,D_{N,r}\}.
\end{equation}


\subsection{Unisolvence and conformity of the finite element}\label{sec:unisolvence}

In this section we will show that we have a properly defined finite element.

\begin{theorem}\label{thm:unisolvence}
  The finite element $\cDS_r(E_N) = \Po_r(E_N)\oplus\Supp_r^\cDS(E_N)$, for $\Supp_r^\cDS(E_N)$
  given by \eqref{eq:supplementSpace}, with nodal DoFs \eqref{eq:nodalFunctionals} is well defined
  (i.e., unisolvent) when $r\ge N-2$. Moreover, a nodal basis is given by the functions defined below in 
  \eqref{eq:nodal-cell}, \eqref{eq:nodal-edge11}, and \eqref{eq:nodal-vertices}.
\end{theorem}

To prove the theorem, we will explicitly construct a basis of shape functions $\varphi_i$ for $\cP$
dual to $\cN$.  Such shape functions are called \emph{nodal basis functions}. For a nodal point
$\x^{\textrm{node}}_j$, they have the property that
$N_j(\varphi_{i})=\varphi_{i}(\x^{\textrm{node}}_j) = \delta_{ij}$, the Kronecker delta.  The
unisolvence property (i.e., that $\cN$ is a basis for the dual space) is then immediate.
Moreover, it follows from the construction that we obtain global $H^1$ conforming elements by just
matching vertex and edge DoFs on the boundaries of the elements; that is, local basis functions
merge together continuously to give a global nodal basis for
$\cDS_r=\cDS_r(\Omega)\subset H^1(\Omega)$. Our construction directly extends that given in
\cite{Arbogast_Tao_Wang_2020x_serendipityMixed} for the case $N=4$.

Before beginning the construction, it is convenient to define
\begin{equation}\label{eq:R-1sided}
	\mathcal{R}_{i,j}(\x) = \tfrac12\big(1 - R_{i,j}(\x)\big),\quad \mathcal{R}_{j,i}(\x) = \tfrac12\big(1 + R_{i,j}(\x)\big), \end{equation}
so that $\mathcal{R}_{k,\ell}$ is $1$ on edge $e_k$, $0$ on $e_\ell$, and arbitrary on the other edges. Let us now set $\lambda_{j,i}=\lambda_{i,j}$ when $i<j$, and define, for any $1\le k,\ell \le N$, $2\le |k-\ell| \le N-2$,
\begin{equation}\label{eq:phi}
  \phi_{k,\ell}(\x) = \Big(\prod_{m\neq k,\ell}\lambda_m\Big)\lambda_{k,\ell}^{r-N+2}\mathcal{R}_{k,\ell}
  \quad\in\cDS_r(E_N).
\end{equation}
These lie in $\Po_r(E)\oplus\Supp_r^{\cDS}(E)$ and satisfy
\begin{align}\label{eq:phi_on_edges}
	\phi_{k,\ell}(\x)=\begin{cases}
		\ \ \ \ \ \ \ \ \ \ \ \quad \quad 0,  &\x \in e_m,\ m \ne k, \\
		\displaystyle\Big(\prod_{m\neq k,\ell}\lambda_m\Big)\lambda_{k,\ell}^{r-N+2}\in \Po_r(e_k), &\x\in e_k.
	\end{cases}
\end{align}
Moreover,
\begin{align}\label{eq:cDS=P+Sextended}
	\cDS_r(E_N) &= \Po_r(E_N) + \spn\{\phi_{k,\ell}:1 \le k,\ell \le N,\ 2\le |k-\ell| \le N-2\}.
\end{align}

\subsubsection{Interior cell nodal basis functions}\label{sec:cell}

For the element $E_N$, we have interior shape functions only when $r\ge N$ (recall Table~\ref{tab:geometricDecomp}). These shape functions are
\begin{equation}\label{eq:nodal-cell-space}
\lambda_1\lambda_2\cdots\lambda_N\Po_{r-N},
\end{equation}
and they vanish on all the edges (i.e., at all edge and vertex nodes).  Let
$\{\phi_{E,i}\}\subset\Po_{r-N}$ be a nodal basis for the cell nodes $\{\x_{E,i}\}$, where
$i=1,2,\ldots,\dim\Po_{r-N}$.  That is, $\phi_{E,i}(\x_{E,j}) = \delta_{ij}$.  Our interior cell nodal basis functions are then
\begin{equation}\label{eq:nodal-cell}
  \varphi_{E,i}(\x) = \frac{\lambda_1(\x)\lambda_2(\x)\cdots\lambda_N(\x)\phi_{E,i}(\x)}
                          {\lambda_1(\x_{E,i})\lambda_2(\x_{E,i})\cdots\lambda_N(\x_{E,i})},
\quad i=1,2,\ldots,\dim\Po_{r-N}.
\end{equation}


\subsubsection{Edge nodal basis functions}\label{sec:edges}

For $\cDS_r(E_N)$, there are $r-1$ edge nodes on each edge. To simplify the notation, we construct
$\varphi_{e,1,1}(x)$, which is $1$ at $x_{e,1,1}$ and vanishes at all other nodal points. The
construction of the other edge nodal basis functions is similar.

For some $\tilde p\in\Po_{r-N+1}(e_1)$ (take $\tilde p = 0$ when $r=N-2$) and for some coefficients
$\beta_j$, let
\begin{equation}\label{eq:nodal-edge-tilde}
  \phi_{e,1,1}(\x) = \Big(\prod_{m\ne 1}\lambda_m(\x)\Big)p(\x) + \sum_{j\ne N,1,2}\beta_j\phi_{1,j}(\x)
  \quad\in\cDS_r(E_N),
\end{equation}
where $p(\x) = \tilde p((\x-\x_{v,N})\cdot\tau_1)$ extends $\tilde p$ to $E_N$ constantly in the normal
direction to $\cL_1$. This function vanishes on all edges but $e_1$.

Denote
\vspace*{-4pt}\begin{equation*}
  t_{e,1,n}=(\x_{e,1,n}-\x_{v,N})\cdot\tau_1
  \quad\text{and}\quad
  \tilde p(t)=\sum_{\ell=0}^{r-N+1}\alpha_\ell\,t^\ell.
\end{equation*}
We require that $\phi_{e,1,1}(\x_{e,1,n})=\delta_{1,n}$ for $n=1,2,\ldots,r-1$, so the $r-1$
coefficients $\{\alpha_\ell,\beta_j\}$ solve the square linear system
\begin{align}\label{eq:nodal-edge-tilde-conditions}
  & \frac{\phi_{e,1,1}(\x_{e,1,n})}{(\lambda_N\lambda_2)(\x_{e,1,n})}
    = \Big(\!\prod_{m\ne N,1,2}\lambda_m(\x_{e,1,n})\Big)\sum_{\ell=0}^{r-N+1}\alpha_\ell\,t_{e,1,n}^\ell\\
  \nonumber
  &\quad
    + \sum_{j\ne N,1,2}\beta_j\Big(\!\prod_{m\ne N,1,2,j}\lambda_m(\x_{e,1,n})\Big)\lambda_{1,j}^{r-N+2}(\x_{e,1,n})
  = \frac{\delta_{1,n}}{(\lambda_N\lambda_2)(\x_{e,1,n})}.
\end{align}

\begin{figure}[ht]
\centerline{
\includegraphics[width=0.304\linewidth]{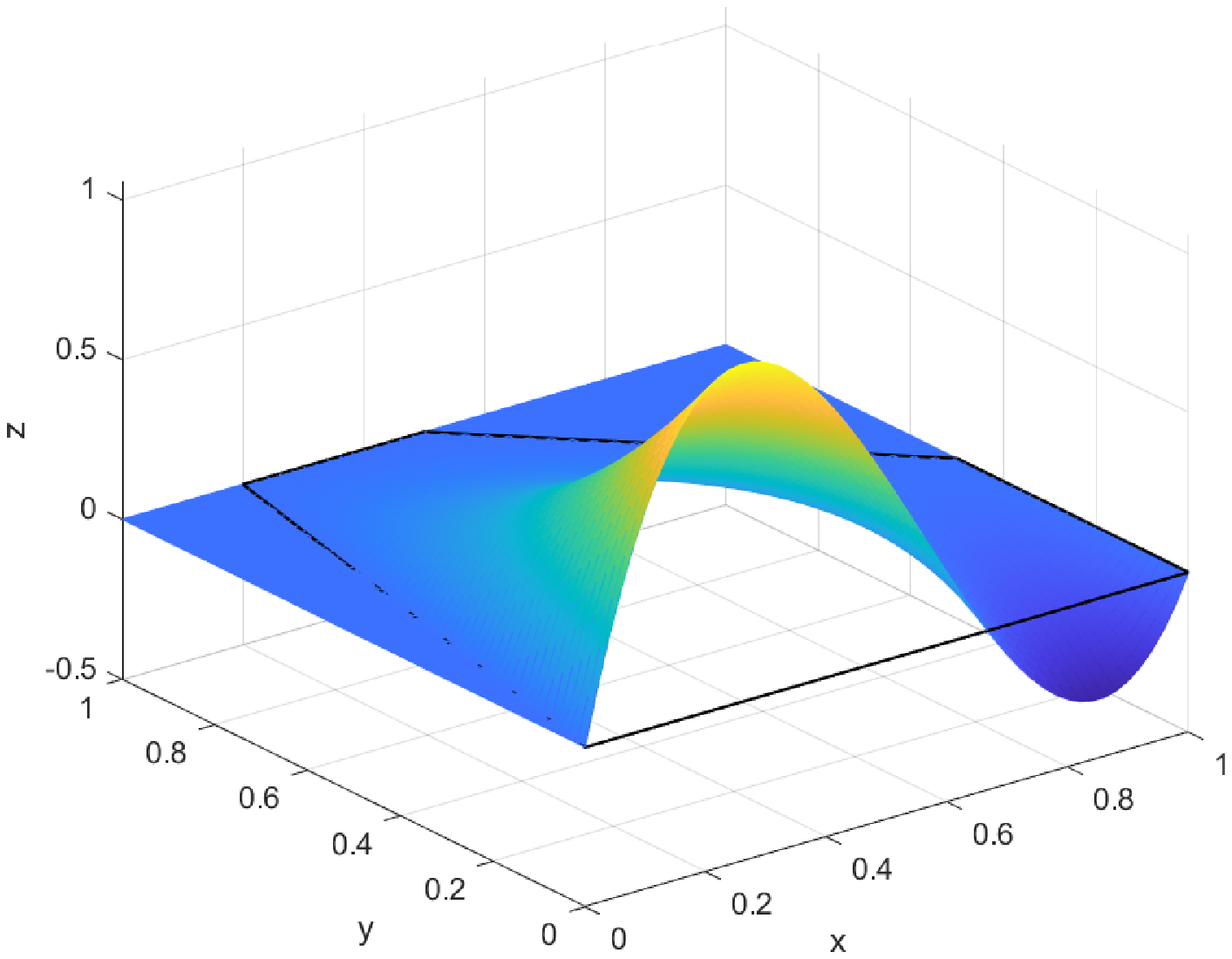}
\qquad\quad\includegraphics[width=0.2\linewidth]{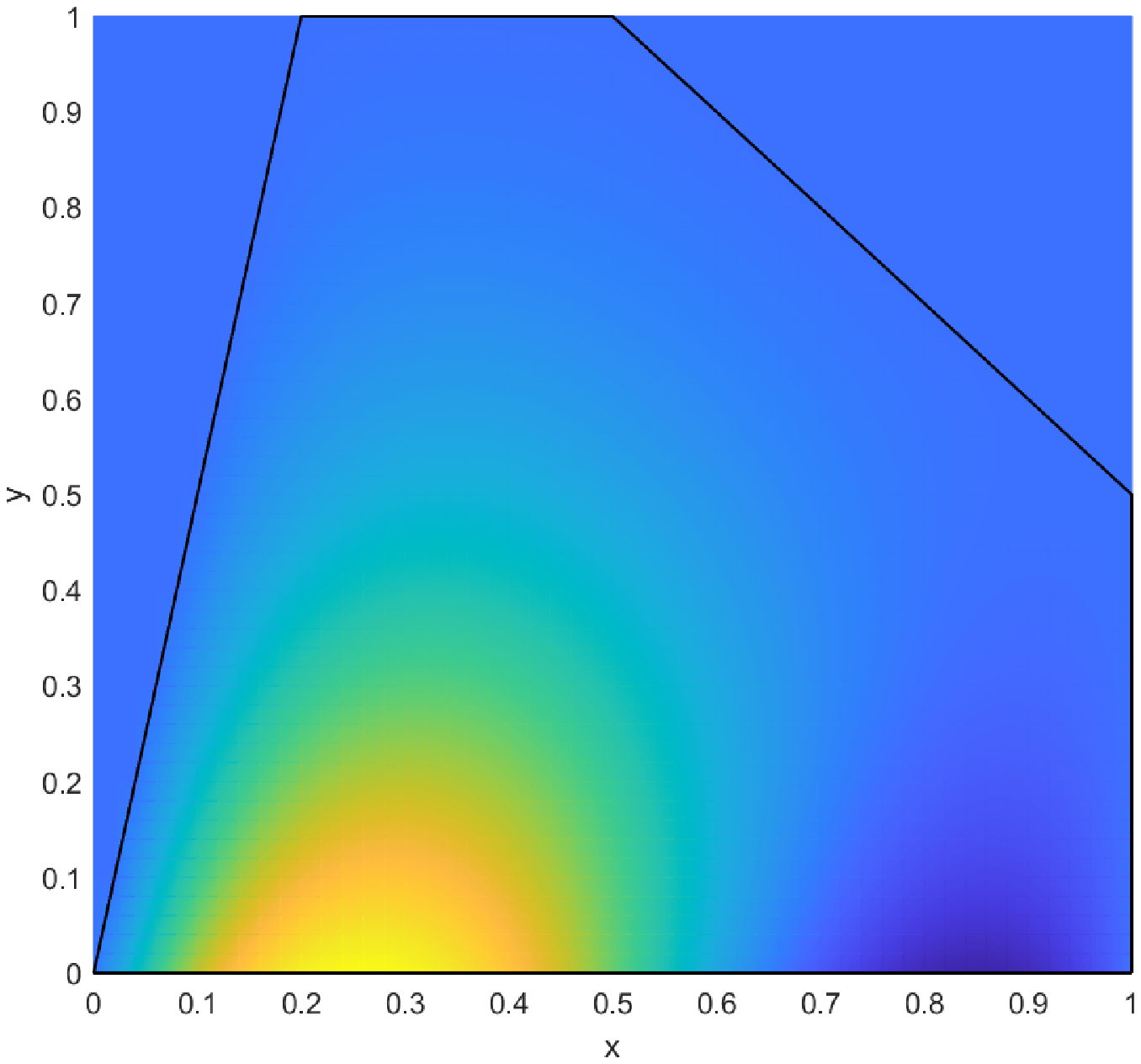}\qquad
}
\caption{Plots of the $r = 3$ basis function for the edge node at $(\frac{1}{3},0)$ of a pentagon. }\label{fig:edge}
\end{figure}

Assume for the moment that the function $\phi_{e,1,1}$ is well defined on $E_N$. It takes the
value $1$ at $\x_{e,1,1}$ and vanishes at all the other vertex and edge nodes, so we define
\begin{equation}\label{eq:nodal-edge11}
  \varphi_{e,1,1}(\x) = \phi_{e,1,1}(\x) - \sum_{k=1}^{\dim\Po_{r-N}(E)}\phi_{e,1,1}(\x_{E,k})\,\varphi_{E,k}(\x).
\end{equation}
The nodal basis functions $\{\varphi_{e,i,j}:i=1,2,\ldots,N,\ j=1,2,\ldots,r-1\}$ for the other edge
nodes are defined similarly. In Figure~\ref{fig:edge}, we show an edge nodal basis function for a
pentagon.

\begin{lemma}\label{lem:solveEdgeLinSyst}
  There exists a unique set of coefficients $\alpha_\ell$, $\ell=0,1,\ldots,r-N+1$, and $\beta_j$,
  $j= 3,4,\ldots,N-1$, solving the $(r-1)\times(r-1)$ linear system \eqref{eq:nodal-edge-tilde-conditions}.
\end{lemma}

\begin{proof} For $t\in\Re$, let $\x(t) = \x_{v,N} + t\,\tau_1$ and define $\tilde q(t)\in\Po_{r-2}(e_1)$ by
\begin{align}\label{eq:nodal-edge-tilde-e-i}
  \tilde q(t) &= \frac{\phi_{e,1,1}(\x(t))}{(\lambda_N\lambda_2)(\x(t))}\\\nonumber
  &= \Big(\!\prod_{m\ne N,1,2}\!\!\lambda_m(\x(t))\Big)\,\tilde p(t)\,
  + \!\!\sum_{j\ne N,1,2}\beta_j\Big(\!\prod_{m\ne N,1,2,j}\!\!\lambda_m(\x(t))\Big)\lambda_{1,j}^{r-N+2}(\x(t)).
\end{align}
We must show that the linear system has a unique solution, which is equivalent to showing that
$\tilde q(t_{e,1,n}) =0$ for all $n=1,2,\ldots,r-1$, then all $\alpha_\ell = 0$ and $\beta_j=0$
($j\ne N,1,2$). Now $\tilde q(t)$ is a polynomial of degree $r-2$, and it vanishes at $r-1$ points,
so it vanishes identically.

Suppose that the lines through $e_1$ and $e_j$ intersect at $\x_{1,j}=\cL_1\cap\cL_j$ for some
$j\ne N,1,2$. Since $\lambda_j(\x_{1,j})=0$, $\tilde q((\x_{1,j}-\x_{v,N})\cdot\tau_1)$ reduces to
\begin{equation*}
0 = \tilde q((\x_{1,j}-\x_{v,N})\cdot\tau_1) 
= \beta_j\Big(\!\prod_{m\ne N,1,2,j}\lambda_m(\x_{1,j})\Big)\lambda_{1,j}^{r-N+2}(\x_{1,j}).
\end{equation*}
But $\lambda_m(\x_{1,j}) \ne 0$ for all $m\ne 1,j$ and $\lambda_{1,j}(\x_{1,j})\ne0$ by our choice of
this linear function, so we conclude that $\beta_j=0$.

We have two cases to consider.  First, if no edge is parallel to $e_1$ (so the intersection points
$\x_{1,j}$ exist for all $j\ne N,1,2$), then all the $\beta_j$ vanish. Second, suppose that the lines
through $e_1$ and $e_{j}$ are parallel for some $j\ne N,1,2$. No other edges can also be parallel,
so we conclude $\beta_k=0$ for all $k\ne j$. Moreover, $\lambda_j|_{e_1}=c>0$ is a strictly positive
constant, and so
\begin{equation*}
  0 = \tilde q(t) = \Big(\!\prod_{m\ne N,1,2,j}\!\!\lambda_m(\x(t))\Big)
  \Big(\sum_{\ell=0}^{r-N+1}\!c\,\alpha_\ell\,t^\ell+\beta_{j}\lambda_{1,j}^{r-N+2}(\x(t))\Big),
\end{equation*}
or
\begin{equation*}
\beta_{j}\lambda_{1,j}^{r-N+2}(\x(t)) = - \sum_{\ell=0}^{r-N+1}\!c\,\alpha_\ell\,t^\ell \in \Po_{r-N+1}(e_1).
\end{equation*}
The zero line of $\lambda_{1,j}$ is transverse to $e_1$ (again by our choice of this linear
function), leading us to conclude that $\lambda_{1,j}^{r-N+2}$ must have strict degree
$r-N+2$. Therefore, again, all the $\beta_{j}=0$.

We have reduced $\tilde q(t)=0$ to a positive function times $\tilde p(t)$, so we must conclude that
$\tilde p(t)=0$.  That is, all the $\alpha_\ell=0$, and the proof is complete.
\end{proof}


\subsubsection{Vertex nodal basis functions}\label{sec:vertexDoFs}

\begin{figure}[ht]
\centerline{
\includegraphics[width=0.304\linewidth]{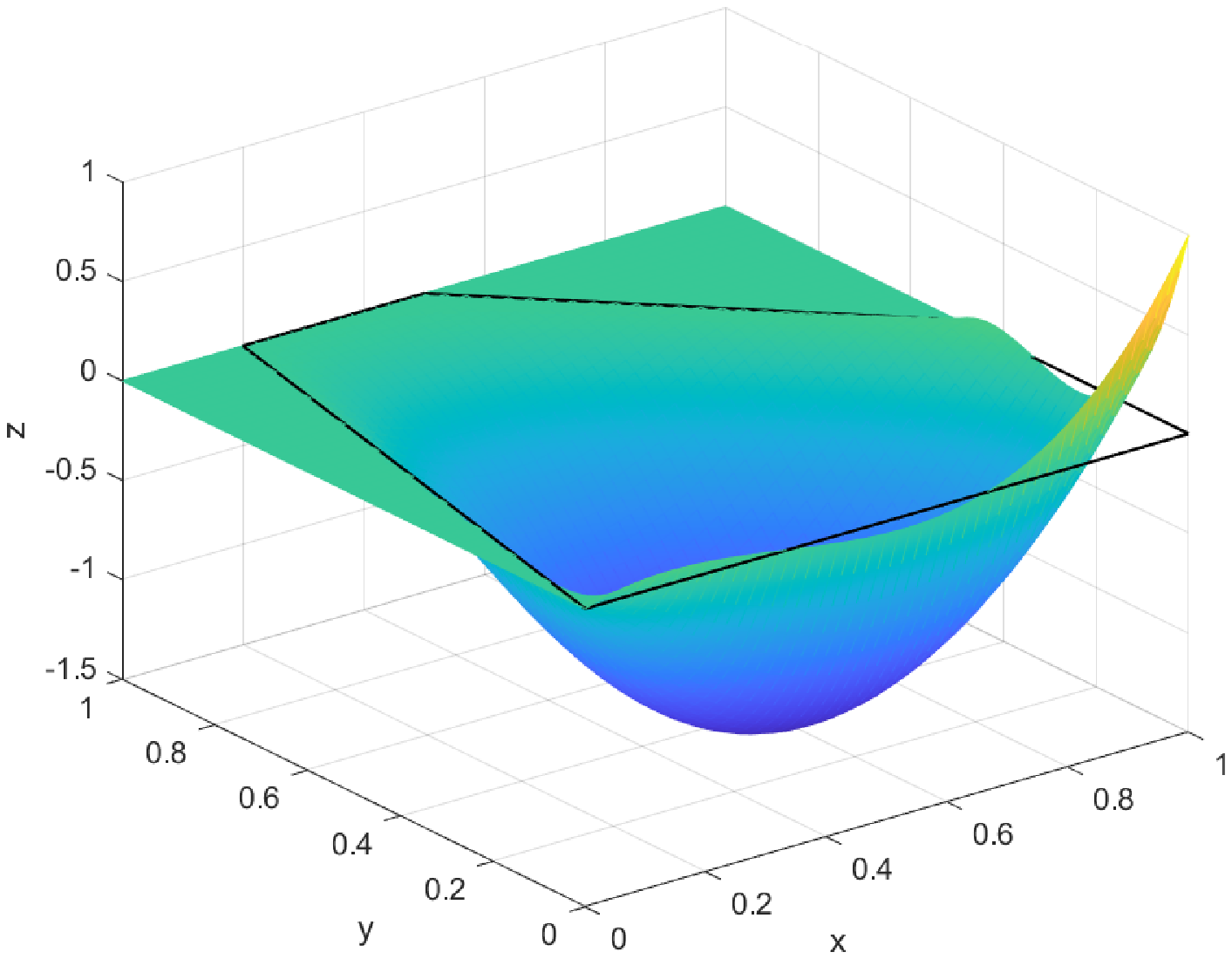}
\qquad\quad\includegraphics[width=0.2\linewidth]{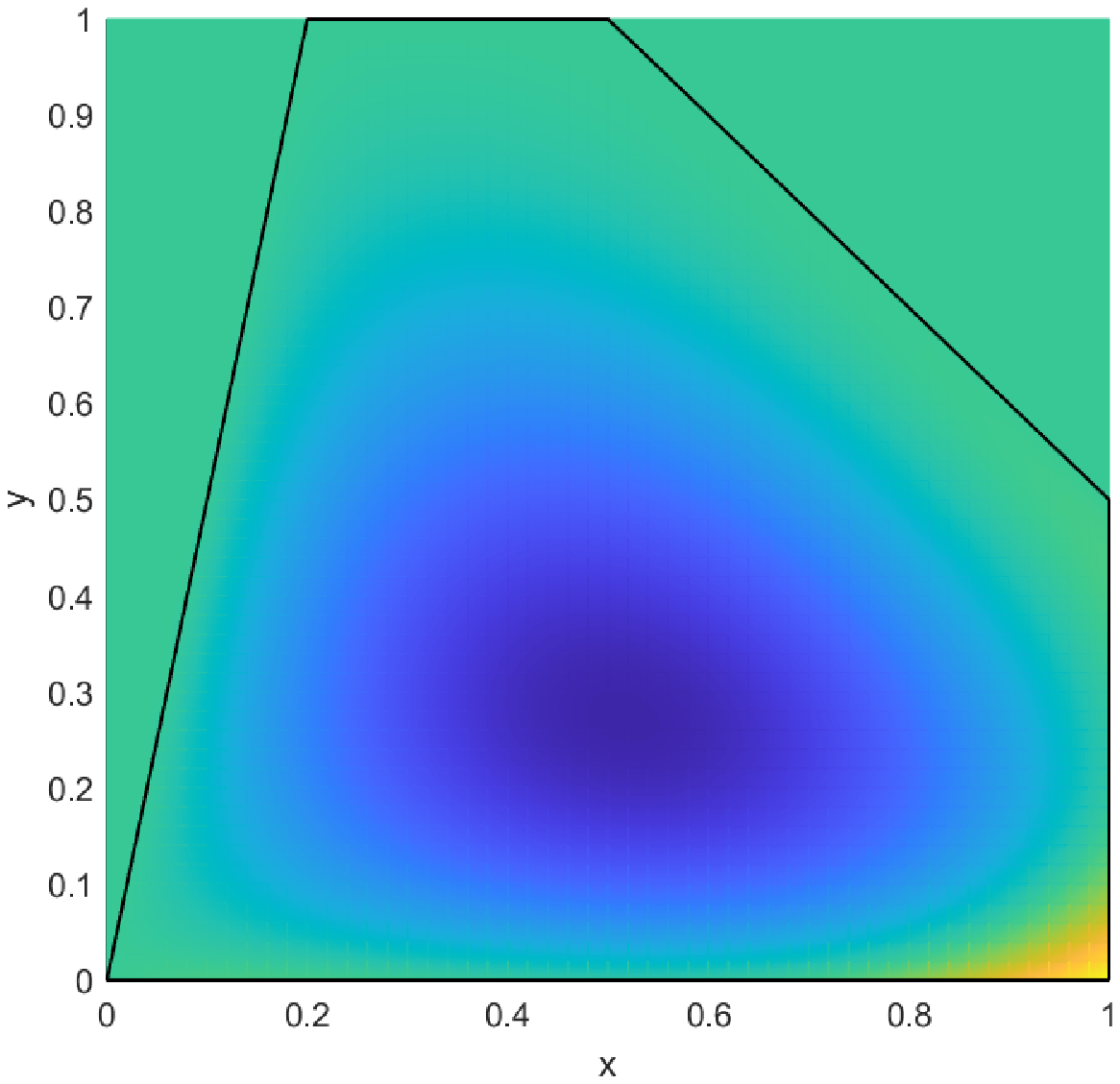}\qquad
}
\caption{Plots of the $r = 3$ basis function for the vertex node at $(1,0)$ of a pentagon. }\label{fig:vertex}
\end{figure}

For the vertices, since $r\ge N-2$, we can define for each $i=1,2,\ldots,N$ the shape functions
\begin{align}\label{eq:nodal-vertices-unnormalized}
  \phi_{v,i}(\x) = \prod_{j\ne i,i+1}\lambda_j(\x)
  - \sum_{k=i}^{i+1}\sum_{\ell=1}^{r-1}\Big(\prod_{j\ne i,i+1}\lambda_j(\x_{e,k,\ell})\Big)\,\varphi_{e,k,\ell}(\x),
\end{align}
wherein we interpret indices modulo $N$. These $N$ functions vanish at all of the edge nodes, and $\phi_{v,i}(\x_{v,j})=0$ if $i\ne j$ and is positive otherwise. The nodal basis functions are then
\begin{align}\label{eq:nodal-vertices}
&\varphi_{v,i}(\x) = \frac{\phi_{v,i}(\x)
       - \sum_{k=1}^{\dim\Po_{r-N}(E)}\phi_{v,i}(\x_{E,k})\,\varphi_{E,k}(\x)}{\phi_{v,i}(\x_{v,i})},
\quad i=1,2,\ldots,N.
\end{align}
A vertex nodal basis function for a pentagon is shown in Figure~\ref{fig:vertex}. This completes the construction of the $D_{N,r}=\dim\Po_r(E_N)+\frac12N(N-3)$ nodal basis functions for $\cDS_r(E_N)$. It also completes the proof of Theorem~\ref{thm:unisolvence}.

\section{Direct serendipity elements when $1\le r<N-2$}\label{sec:DSlow}

There are vertex and possibly edge nodes, but no interior nodes,
when $1\le r<N$. The total number of DoFs needed for $E_N$ is then simply
\begin{equation}\label{eq:totalNumber_small_r}
D_{N,r} = N + N(r-1) = Nr \ge \dim\Po_r(E_N).
\end{equation}
Our strategy is to define the space as a subset of a higher order direct serendipity space; that is,
for some index $s$ such that $r<s<N$, we define
\begin{align}\label{eq:cDS_small_r}
  \cDS_r^{(s)}(E_N)
  = \big\{\varphi\in\cDS_s(E_N) : \varphi|_e\in\Po_r(e)\text{ for all edges $e$ of }E_N\big\}.
\end{align}

\begin{theorem}\label{thm:DS_low_r}
  The finite element \eqref{eq:cDS_small_r} with nodal DoFs \eqref{eq:nodalFunctionals}  is well defined
  (i.e., unisolvent) when $r<N-2$ and $r<s<N$. Moreover,  
  \begin{align}\label{eq:cDS=P+S_small_r}
    \cDS_r^{(s)}(E_N) = \Po_r(E_N)\oplus\Supp_r^\cDS(E_N)
  \end{align}
  for some supplemental space of functions $\Supp_r^\cDS(E_N)$, and a nodal basis is given by the
  functions listed in \eqref{eq:ds_sm_r} and defined as in \eqref{eq:varphi_r_e11} and
  \eqref{eq:varphi_r_v1}.
\end{theorem}

As a practical matter, one should take $s=N-2$. It is obvious that
$\Po_r(E_N)\subset\cDS_r^{(s)}(E_N)$, since $\Po_r(E_N)\subset\Po_s(E_N)\subset\cDS_s(E_N)$
restricts to $\d E_N$ as required. That is, $\cDS_r^{(s)}(E_N)$ has the form
\eqref{eq:cDS=P+S_small_r}.  We prove the rest of the theorem in the next section by constructing a
nodal basis.


\subsection{Construction of the nodal basis functions when $r<N-2$}\label{sec:universal_small_r}

We construct nodal basis functions for $\cDS_r(E_N)$ from
$\cDS_s(E_N)$ for any $r<s<N$. To make the notation clear as to which order ($r$ or $s$) a quantity refers to, we will use a superscript within parentheses.  For example, edge node $\x_{e,1}$ will be referred to as $\x_{e,1,1}^{(r)}$ if it is the node in $\cDS_r(E)$, and $\x_{e,1,1}^{(s)}$ if it is the node in $\cDS_s(E)$ (these two nodes are not at the same position).

We first note that for each $j=0,1,\ldots,r$, there exists a unique
$\tilde p^{(r)}_{j}(t)\in\Po_r([0,1])$ interpolating $r+1$ points as
\begin{equation}\label{eq:e11interpPoly}
	\tilde p^{(r)}_{j}(k/r) = \delta_{j,k},\quad \forall k=0,1,\ldots,r.
      \end{equation}

A basis function for edge node $\x_{e,i,j}^{(r)}$, $i=1,2,\ldots,N$ and $j=1,2,\ldots,r-1$, is then
\begin{equation}\label{eq:varphi_r_e11}
	\varphi^{(r,s)}_{e,i,j}(\x) = \sum_{k=1}^{s-1}\tilde p^{(r)}_{j}(k/s)\,\varphi^{(s)}_{e,i,k}(\x)\in\cDS_s(E_N),
\end{equation}
which vanishes on all the edges except for $e_i$.  Restricted to $e_i$, it is nominally a polynomial
of degree~$s$. However, it agrees with $\tilde p^{(r)}_{j}$ at $s+1>r+1$ points, so it is in fact a
polynomial of degree $r$ on $e_i$.  In consequence, $\varphi^{(r,s)}_{e,i,j}\in\cDS_r^{(s)}(E_N)$,
and it vanishes at all nodes of $\cDS_r^{(s)}(E)$ except $\x_{e,i,j}^{(r)}$, where it is one (i.e.,
it is a nodal basis function).

For a vertex node $\x_{v,i}=\x_{v,i}^{(r)}=\x_{v,i}^{(s)}$, we define
\begin{align}\label{eq:varphi_r_v1}
	\varphi^{(r,s)}_{v,i}(\x)
	&= \varphi^{(s)}_{v,i}(\x) + \sum_{j=1}^{s-1}\tilde p^{(r)}_{r}(j/s)\,\varphi^{(s)}_{e,i,j}(\x)
         + \sum_{j=1}^{s-1}\tilde p^{(r)}_{0}(j/s)\,\varphi^{(s)}_{e,i+1,j}(\x)\\\nonumber
	&\in\cDS_s(E_N),
\end{align}
which vanishes on all the edges except $e_i$ and $e_{i+1}$. As before, we conclude that it is a
polynomial of degree $r$ on edges $e_i$ and $e_{i+1}$, and so
$\varphi^{(r,s)}_{v,i}(\x)\in\cDS_r^{(s)}(E_N)$. Moreover, it is the nodal basis function for
$\x_{v,i}$, since it vanishes at all edge nodes $\x_{e,k,j}^{(r)}$ of $e_k$, $k=i,i+1$, and
$\varphi^{(r,s)}_{v,i}(\x_{v,i})=1$.

Finally, since there are no interior cell DoFs, we conclude that
\begin{align}\label{eq:ds_sm_r}
  \cDS_r^{(s)}(E_N)  &= \spn\big\{\{\varphi^{(r,s)}_{v,i}: i=1,2,\ldots,N\}\\\nonumber
  &\qquad\qquad\cup \{\varphi^{(r,s)}_{e,i,j}: i=1,2,\ldots,N,\ j = 1,2,\ldots,r-1 \}\big\},
\end{align}
which indeed has dimension $Nr$.  This completes the proof of Theorem~\ref{thm:DS_low_r}.
  

\subsection{A second construction identifying the supplemental function space}\label{sec:whole_space_small_r}

From either the definition \eqref{eq:cDS_small_r} or from the nodal basis \eqref{eq:ds_sm_r}, it is
difficult to determine the supplemental space $\Supp_r^\cDS(E_N)$ in \eqref{eq:cDS=P+S_small_r}. In
this section, we give an explicit construction $\Supp_r^\cDS(E_N)$. In practice, the supplemental
space is not needed to implement $\cDS_r(E_N)$ (one would simply use \eqref{eq:ds_sm_r}); however,
as we will see later, it could be used to implement mixed finite elements.

It will be convenient in this section to use a notation that unifies edge and vertex
nodes. For each edge index $i=1,2,\ldots,N$ and $j=0,1,\ldots,r$, let
\begin{equation}
  \x_{n,i,j} = \begin{cases}
    \x_{v,i-1}&\text{if }j=0\text{ (where $i-1$ is interpreted as $N$ when $i=1$)},\\
    \x_{v,i}&\text{if }j=r,\\
    \x_{e,i,j}&\text{if }j=1,2,\ldots,r-1.
  \end{cases}
\end{equation}
We caution that the vertices are represented twice in this indexing convention. Let the full set of
nodal points be denoted
\begin{align*}
  \cA &= \{\x_{v,i}, \x_{e,i,j} : i=1,2,\ldots,N,\ j=1,2,\ldots,r-1\}\\
  &= \{\x_{n,i,j}:i=1,2,\ldots,N,\ j=1,2,\ldots,r\}.
\end{align*}
We will divide this set into two disjoint subsets $\cA_\Po$ and $\cA_\Supp=\cA\setminus\cA_\Po$.

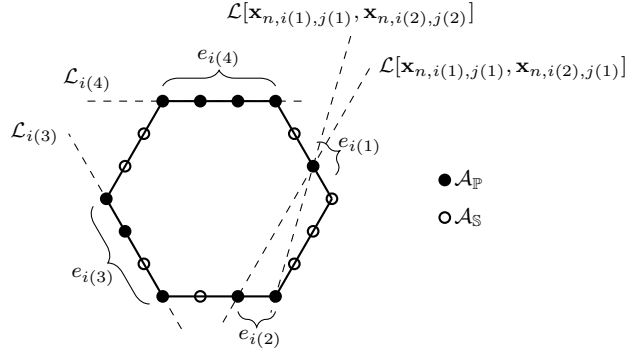
\begin{figure}\label{fig:SuppFunc_smallr}
	\parbox{\textwidth}{\centering\begin{tikzpicture}[xscale=0.25,yscale=0.25]
		\draw [thick, -] (-3,0) -- (3,0);
		\draw [thick, -] (3,0) -- (6,5.196);
		\draw [thick, -] (6,5.196) -- (3,10.392);
		\draw [thick, -] (3,10.392) -- (-3,10.392);
		\draw [thick, -] (-3,10.392) -- (-6,5.196);
		\draw [thick, -] (-6,5.196) -- (-3,0);
		%
		%
		%
		\draw [thick] (-1,0) circle (8pt);
		\draw [thick] (4,1.732) circle (8pt);
		\draw [thick] (5,3.464) circle (8pt);
		\draw [thick] (6,5.196) circle (8pt);
		\draw [thick] (4,8.660) circle (8pt);
		\draw [thick] (-4,1.732) circle (8pt);
		\draw [thick] (-5,6.928) circle (8pt);
		\draw [thick] (-4,8.660) circle (8pt);
		\filldraw [thick] (5,6.928) circle (8pt); 
		\draw [decorate, decoration={brace,amplitude=6pt,mirror}] (6.193,6.462) -- (5.193,8.194);
		\draw (7.5,7) node [above] {\footnotesize $e_{i(1)}$};
		\filldraw [thick] (1,0) circle (8pt); 
		\filldraw [thick] (3,0) circle (8pt);
		\draw [decorate, decoration={brace,amplitude=6pt,mirror}] (1,-0.8) -- (3,-0.8);
		\draw (2.2,-1.3) node [below] {\footnotesize $e_{i(2)}$};
		\filldraw [thick] (-3,0) circle (8pt);
		\filldraw [thick] (-5,3.464) circle (8pt);
		\filldraw [thick] (-6,5.196) circle (8pt);
		\draw [decorate, decoration={brace,amplitude=6pt}] (-3.693,-0.4) -- (-6.693,4.796);
		\draw (-6.693,2) node [below] {\footnotesize $e_{i(3)}$};
		\filldraw [thick] (-3,10.392) circle (8pt);
		\filldraw [thick] (-1,10.392) circle (8pt);
		\filldraw [thick] (1,10.392) circle (8pt);
		\filldraw [thick] (3,10.392) circle (8pt);
		\draw [decorate, decoration={brace,amplitude=6pt}] (-3,11.192) -- (3,11.192);
		\draw (0,11.7) node [above] {\footnotesize $e_{i(4)}$};
		\draw [dashed] (-7,10.392) node[above]{\footnotesize$\cL_{i(4)}$} -- (4.5, 10.392) ;
		\draw [dashed] (-2,-1.732)  -- (-8, 8.660) node[left]{\footnotesize$\cL_{i(3)}$};
		\draw [dashed] (8,12.134) node[right]{\footnotesize$\cL[\x_{n,i(1),j(1)},\x_{n,i(2),j(1)}]$} -- (0, -1.732);
		\draw [dashed] (7,13.856) node[above]{\footnotesize$\cL[\x_{n,i(1),j(1)},\x_{n,i(2),j(2)}]$} -- (2.75, -0.866);
		\filldraw [thick] (12,6.196) circle (8pt) node[right]{\footnotesize$\cA_\Po$};
		\draw [thick] (12,4.196) circle (8pt) node[right]{\footnotesize$\cA_\Supp$};
		\end{tikzpicture}}
	\caption{A choice of nodes $\cA_\Po$ and $\cA_\Supp$ for $N=6$, $r=3$. And the dashed lines show the choices of zero lines for the construction of $\phi_{n,i(2),j(\ell)}$, $\ell=1,2$.}
	\label{Fig:SuppFunc_smallr}
\end{figure}

The subset of nodes $\cA_\Po$ is chosen iteratively as follows, and as depicted in Figure~\ref{fig:SuppFunc_smallr}. For each $k=r+1,\ldots,2,1$ in descending order, first select a distinct edge $e_{i(k)}$ with index $i(k)\in\{1,2,\ldots,N\}$. At this stage, there are at least $N-r+k-1>0$ edges left to choose from, since $N-r>2$ and $k\ge1$. Second, select $k$ distinct nodes $\x_{n,i(k),j(\ell)}$ on this chosen edge, with the indices $j(\ell)\in\{0,1,\ldots,r\}$ and $\ell=1,2,\ldots,k$. The only restriction is that one may not choose a vertex node that lies on any of the previously chosen edges. Since there are 2 vertex nodes and $r-1$ edge nodes, one can always meet this restriction. As a simple example, one can choose edges $i(k)=k$ and take only edge nodes, except for $\x_{v,r}$ and $\x_{v,r+1}$ on $e_{r+1}$ and $\x_{v,r-1}$ on $e_{r}$.

The total number of nodes in $\cA_\Po$ is
\begin{equation}\label{eq:dim_Po_r}
  \sum_{k=1}^{r+1}k = \frac12(r+2)(r+1) = \dim\Po_r(E_N).
\end{equation}
The total number of unselected nodes $\cA_\Supp=\cA\setminus\cA_\Po$ is the same as the dimension of
$\Supp_r^\cDS(E_N)$. For each node $\x_{n,i,j}\in\cA_\Supp$, we construct
$\varphi_{n,i,j}=\varphi_{n,i,j}^{(r,s)}$, the supplemental function associated to  $\x_{n,i,j}$ as
in the previous section.  The supplemental space is then
\begin{equation}\label{eq:supp_small_r}
  \Supp_r^\cDS(E_N) = \spn\{\varphi_{n,i,j} : \x_{n,i,j}\in\cA_\Supp\}\subset\cDS_r^{(s)}(E_N),
\end{equation}
and it has the correct dimension.  These basis functions are nodal, by construction.

To verify that \eqref{eq:supp_small_r} is indeed the supplemental space, we finish the construction of
the nodal basis (i.e., for nodal points in $\cA_\Po$) by including additional functions only from
$\Po_r(E_N)$. We do this iteratively for each $k=1,2,\ldots,r+1$ in ascending order as follows.  For
$k=1$, we construct the nodal basis function for $\x_{n,i(1),j(1)}$ by first defining
\begin{equation}\label{eq:n_i1_j1_pre}
  \phi_{n,i(1),j(1)}(\x) = \prod_{m=2}^{r+1}\frac{\lambda_{i(m)}(\x)}{\lambda_{i(m)}(\x_{n,i(1),j(1)})}\ \in\Po_r,
\end{equation}
which vanishes at all the nodes of $\cA_\Po$ except $\x_{n,i(1),j(1)}$, where it is one. By the choice of edges, the denominator does not vanish. Then
\begin{equation}\label{eq:n_i1_j1}
  \varphi_{n,i(1),j(1)}(\x)
  = \phi_{n,i(1),j(1)}(\x) - \sum_{\x_{n,i,j}\in\cA_\Supp}\phi_{n,i(1),j(1)}(\x_{n,i,j})\,\varphi_{n,i,j}(\x),
\end{equation}
and this is indeed our nodal basis function for the node $\x_{n,i(1),j(1)}$.

For $k=2$, we need to construct the nodal basis functions for the two points on the edge $e_{i(2)}$. Note that we have one more point compared to the previous step, but we also have one fewer edge to deal with, since we now have $\varphi_{n,i(1),j(1)}$. Therefore we can construct for each $\ell=1,2$,
\begin{equation*}
  \phi_{n,i(2),j(\ell)}(\x)
  = \frac{\lambda[\x_{n,i(1),j(1)},\x_{n,i(2),j(\ell^*)}](\x)}{\lambda[\x_{n,i(1),j(1)},\x_{n,i(2),j(\ell^*)}](\x_{n,i(2),j(\ell)})}
                           \prod_{m=3}^{r+1}\frac{\lambda_{i(m)}(\x)}{\lambda_{i(m)}(\x_{n,i(2),j(\ell)})}\ \in\Po_r,
\end{equation*}
where $\ell^*=2,1$ is the other index. For each $\ell=1,2$, the function vanishes at all the nodes of $\cA_\Po$ except $\x_{n,i(2),j(\ell)}$, where it is one. Then let
\begin{align*}
  \varphi_{n,i(2),j(\ell)}(\x)
  &= \phi_{n,i(2),j(\ell)}(\x) - \!\sum_{\x_{n,i,j}\in\cA_\Supp}\!\phi_{n,i(2),j(\ell)}(\x_{n,i,j})\,\varphi_{n,i,j}(\x),
\end{align*}
which give our two desired nodal basis functions on $e_{i(2)}$.

Perhaps the general construction is clear.  For $k=1,2,\ldots,r+1$, first define for each $\ell=1,2,\ldots,k$,
\begin{align*}
  &\phi_{n,i(k),j(\ell)}(\x)\\
  &\quad= \prod_{m=1,m\ne\ell}^k\frac{\lambda[\x_{n,i(1),j(1)},\x_{n,i(k),j(m)}](\x)}
                                    {\lambda[\x_{n,i(1),j(1)},\x_{n,i(k),j(m)}](\x_{n,i(k),j(\ell)})}
     \prod_{m=k+1}^{r+1}\frac{\lambda_{i(m)}(\x)}{\lambda_{i(m)}(\x_{n,i(k),j(\ell)})}\ \in\Po_r,
\end{align*}
and then set
\begin{align*}
  \varphi_{n,i(k),j(\ell)}(\x)
  &= \phi_{n,i(k),j(\ell)}(\x) - \!\sum_{\x_{n,i,j}\in\cA_\Supp}\!\phi_{n,i(k),j(\ell)}(\x_{n,i,j})\,\varphi_{n,i,j}(\x)\\
   &\qquad - \sum_{m=2}^{k-1}\sum_{l=1}^m\phi_{n,i(k),j(\ell)}(\x_{n,i(m),j(l)})\,\varphi_{n,i(m),j(l)}(\x).
\end{align*}
This completes the identification of $\cDS_r^{(s)}(E_N)$ as $\Po_r(E_N)\oplus\Supp_r^\cDS(E_N)$ for the
supplemental space defined by \eqref{eq:supp_small_r}.


\section{Approximation properties of $\cDS_r$}\label{sec:approxDS}

To obtain global approximation properties, we need to assume that the mesh is uniformly shape
regular in some sense.  We take the definition due to Girault and Raviart
\cite[pp.~104--105]{Girault_Raviart_1986}.

\begin{definition}\label{defn:shape_regular}
For any $E_N\in\cT_h$,
denote by $T_i$, $i=1,2,\ldots,N(N-1)(N-2)/6$, the sub-triangle of $E_N$ with vertices being three of the $N$ vertices
of $E_N$. Define the parameters
\begin{align}
\label{eq:hE}
h_{E_N} &= \text{diameter of }E_N, \\
\label{eq:rhoE}
\rho_{E_N} &= 2\,\min_{1\leq i\leq N(N-1)(N-2)/6}\{ \text{diameter of largest circle inscribed in }T_i \}.
\end{align}
A collection of meshes $\{\cT_h\}_{h>0}$ is \emph{uniformly shape regular} if there exists a
\emph{shape regularity parameter} $\sigma_*>0$, independent of $\cT_h$ and $h>0$, such that the
ratio
\begin{equation}
  \frac{\rho_{E_N}}{h_{E_N}} \geq \sigma_*>0\quad\text{for all }E_N\in\cT_h.
\end{equation}
\end{definition}

A shape regular mesh has the property that every element can take on vertices only in a compact set
of possible values (up to translation and rotation). It also has a bound on the number of elements
that can share a single vertex.  We need the following hypothesis on the construction of the finite elements.

\begin{assumption}\label{assumption:niceDS}
  For every $E_N\in\cT_h$, suppose that the basis functions of $\cDS_{r}(E_N)$ are constructed using
  $\lambda_{i,j}$ such that the zero set $\cL_{i,j}$ intersects $e_i$ and $e_j$. Moreover, assume
  that ${R}_{i,j}$ are uniformly differentiable functions of the vertices of $E_N$ up to order~$m\leq r+1$.
\end{assumption}

\begin{theorem}\label{thm:approxDS}
  Let $\cT_h$ be uniformly shape regular with shape regularity parameter $\sigma_*$ and let
  Assumption~\ref{assumption:niceDS} hold. Let $1\le p\le\infty$ and $\ell>1/p$ (or $\ell\ge1$ if
  $p=1$). Then for $r\ge1$, there exists a constant $C=C(r, \sigma_*)>0$, independent of
  $h=\max_{E_N\in\cT_h}h_{E_N}$, such that for all functions $v\in W^{\ell,p}(\Omega)$,
  \begin{align}
    \inf_{v_h\in\cDS_r(\Omega)}\|v - v_h\|_{W^{m,p}(\Omega)} \leq C\,h^{\ell-m}\,\|v\|_{W^{\ell,p}(\Omega)},
    \quad 0\leq m \leq \ell \leq r+1.
  \end{align}
\end{theorem}

The proof follows closely that given in \cite{Arbogast_Tao_Wang_2020x_serendipityMixed} for the
quadrilateral case and so is omitted here except for discussion of one important issue. The proof
uses a continuous dependence argument, relying on the fact that the set of vertices lies in a
compact set as well as Assumption~\ref{assumption:niceDS}, which ensures that the construction of
the finite elements on $E_N$ depends continuously on its vertices.  The issue that arises when
dealing with polygons is settling on a suitable reference configuration, from which the true element
of the mesh is a continuous and compact perturbation.

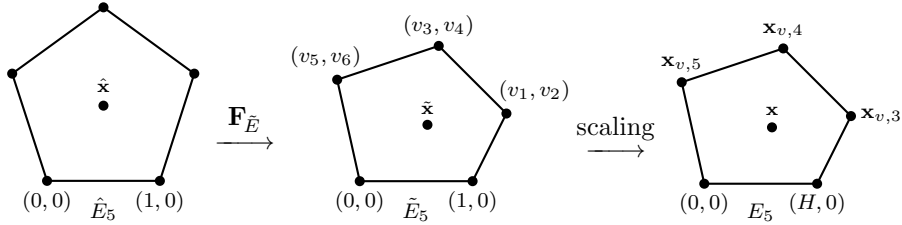
\begin{figure}[ht]
\centering{\setlength\unitlength{3.5pt}\begin{tikzpicture}[xscale=1.5,yscale=1.5]
\draw [thick, -] (0,0) -- (1,0);
\draw [thick, -] (1,0) -- (1.305,0.951);
\draw [thick, -] (1.305,0.951) -- (0.5,1.539);
\draw [thick, -] (0.5,1.539) -- (-0.309,0.951);

\draw [thick, -] (-0.309,0.951) -- (0,0);

\filldraw [thick] (0,0) circle (1pt) node[below]{\footnotesize$(0,0)$};
\node [thick] at (0.5,-0.25) {\footnotesize$\hat E_5$};

\filldraw [thick] (1,0) circle (1pt) node[below]{\footnotesize$(1,0)$};

\filldraw [thick] (1.305,0.951) circle (1pt);

\filldraw [thick] (0.5,1.539) circle (1pt);

\filldraw [thick] (-0.309,0.951) circle (1pt);

\filldraw [thick] (0.5,0.667) circle (1pt);

\node [thick] at (0.5,0.85) {\scriptsize $\hat \x$};
\end{tikzpicture}
\ 
\begin{picture}(6,21)(-3,-11)\small
\thicklines
\put(0,0){\makebox(0,0){$\overset{\text{\normalsize$\F_{\!\tilde E}$}}{-\!\!\!-\!\!\!\longrightarrow}$}}
\end{picture}
\begin{tikzpicture}[xscale=1.5,yscale=1.5]
\draw [thick, -] (0,0) -- (1,0);
\draw [thick, -] (1,0) -- (1.3,0.6);
\draw [thick, -] (1.3,0.6) -- (0.7,1.2);
\draw [thick, -] (0.7,1.2) -- (-0.2,0.9);
\draw [thick, -] (-0.2,0.9) -- (0,0);

\filldraw [thick] (0,0) circle (1pt) node[below]{\footnotesize$(0,0)$};
\node [thick] at (0.5,-0.25) {\footnotesize$\tilde E_5$};

\filldraw [thick] (1,0) circle (1pt) node[below]{\footnotesize$(1,0)$};

\filldraw [thick] (1.3,0.6) circle (1pt);
\node [thick] at (1.55,0.8) {\footnotesize$(v_1,v_2)$};

\filldraw [thick] (0.7,1.2) circle (1pt) node[above]{\footnotesize$(v_3,v_4)$};

\filldraw [thick] (-0.2,0.9) circle (1pt);
\node [thick] at (-0.3,1.1) {\footnotesize$(v_5,v_6)$};

\filldraw [thick] (0.6,0.5) circle (1pt);

\node [thick] at (0.6,0.65) {\scriptsize $\tilde \x$};

\end{tikzpicture}
\begin{picture}(6,20)(-3,-10)\small
\thicklines
\put(0,0){\makebox(0,0){$\overset{\text{\normalsize scaling}}{-\!\!\!-\!\!\!\longrightarrow}$}}
\end{picture}
\begin{tikzpicture}[xscale=1.5,yscale=1.5]
\draw [thick, -] (0,0) -- (1,0);
\draw [thick, -] (1,0) -- (1.3,0.6);
\draw [thick, -] (1.3,0.6) -- (0.7,1.2);
\draw [thick, -] (0.7,1.2) -- (-0.2,0.9);
\draw [thick, -] (-0.2,0.9) -- (0,0);

\filldraw [thick] (0,0) circle (1pt) node[below]{\footnotesize$(0,0)$};
\node [thick] at (0.5,-0.25) {\footnotesize$E_5$};

\filldraw [thick] (1,0) circle (1pt) node[below]{\footnotesize$(H,0)$};

\filldraw [thick] (1.3,0.6) circle (1pt) node [right] {\footnotesize$\x_{v,3}$};

\filldraw [thick] (0.7,1.2) circle (1pt) node [above] {\footnotesize$\x_{v,4}$};

\filldraw [thick] (-0.2,0.9) circle (1pt) node [above] {\footnotesize$\x_{v,5}$};

\filldraw [thick] (0.6,0.5) circle (1pt);

\node [thick] at (0.6,0.65) {\scriptsize $\x$};

\end{tikzpicture}}
\caption{An element $E_5\in\cT_h$ is shown on the right-hand side in its translated and rotated
  local coordinates.  It is the image of a regular reference polygon $\hat E_5$ on the left-hand
  side.  The map is decomposed into one that changes the geometry but not the size
  $\F_{\!\tilde E}:\hat E_5\to\tilde E_5$, and a scaling map
  $\tilde\x\mapsto H\tilde\x$.\label{fig:localRefE}}
\end{figure}

The main argument is illustrated in Figure~\ref{fig:localRefE} for a pentagonal element
$E_N=E_5\in\cT_h$ for which, after translation and rotation, $\x_{v,1}=(0,0)$ and $\x_{v,2}=(H,0)$.
The reference domain is a regular polygon (equilateral and equiangular) $\hat E_N$ with two fixed
vertices $\hat\x_{v,1}=(0,0)$ and $\hat\x_{v,2}=(1,0)$. We need a bijective and smooth map
$\F_{\!\tilde E_N}:\hat E_N\to\tilde E_N=E_N/H$ with $(\ell,0)$ being mapped to $(\ell,0)$,
$\ell=0,1$. In the case of a quadrilateral, one uses a bilinear map.  For a polygon, it is probably
clear to the reader that such a map $\F_{\!\tilde E_N}:\hat E_N\to\tilde E_N$ exists. To be
rigorous, however, we construct $\F_{\!\tilde E_N}$ using smooth barycentric coordinates
$\{\hat\varphi_{v,i}:i=1,\ldots,N\}$ on $E_N$ \cite{Floater_Hormann_Kos_2006}.  The map is then
\begin{equation}\label{eq:Fmap}
  \F_{\!\tilde E_N}(\hat\x) = \sum_{i=1}^N \frac1H \x_{v,i}\,\hat\varphi_{v,i}(\hat\x).
\end{equation}


\section{The de Rham complex and mixed finite elements}\label{sec:deRham}

The de Rham complex of interest here is
\begin{equation}\label{eq:deRham}
\Re \hooklongrightarrow H^1 \overset{\Curl\,}{\longlongrightarrow}
H(\Div) \overset{\Div\,}{\longlongrightarrow} L^2 \longrightarrow 0,
\end{equation}
where the curl (or rot) of a scalar function $\phi(\x)=\phi(x_1,x_2)$ is
$\Curl\,\phi = \bigg(\dfrac{\d\phi}{\d x_2},-\dfrac{\d\phi}{\d x_1}\bigg)$.  From left to right, the
image of one linear map is the kernel of the next.


\subsection{Direct mixed finite elements on polygons}\label{sec:directMixedFE}

For $r=0$ and $s=0$, as well as for each $r=1,2\ldots$ and $s=r-1,r$, there are important discrete
analogues of the de Rham complex involving the direct serendipity spaces and mixed finite element
spaces, denoted $\V_r^{s}(E_N)$, namely
\begin{equation}\label{eq:deRhamV}
  \Re \hooklongrightarrow \cDS_{r+1}(E_N) \overset{\Curl\,}{\longlongrightarrow}
  \V_r^{s}(E_N) \overset{\Div\,}{\longlongrightarrow} \Po_{s}(E_N) \longrightarrow 0.
\end{equation}
On triangular and rectangular elements when $r\geq1$, it is known that the classic serendipity space
$\cS_{r+1}$ (in place of $\cDS_{r+1}$ above) is the precursor of the Brezzi-Douglas-Marini mixed
finite element space BDM$_r$ \cite{BDM_1985,Arnold_Awanou_2011,Arnold_Awanou_2014} (in place of
$\V_r^{r-1}$ above). It is also known that on quadrilateral elements, the direct serendipity space is
the precursor of the direct mixed spaces \cite{Arbogast_Tao_Wang_2020x_serendipityMixed}. The
families of mixed finite elements on $E_N$, $N>4$, are new.

To dissect the properties of these new elements, we note two facts. First, the divergence operator
takes $\x\Po_s$ one-to-one and onto $\Po_s$. Second, the well-known Helmholtz-like decomposition
holds~\cite{Arbogast_Correa_2016}
\begin{equation}
\Po_r^2 = \Curl\,\Po_{r+1}\oplus\x\Po_{r-1}.
\end{equation}
From \eqref{eq:deRham}, we have a reduced ($s=r-1\ge0$) and full ($s=r$) $H(\Div)$-approximating
mixed finite element space ($\cP$ in Definition~\ref{defn:ciarlet}) defined directly on a polygon
$E_N$ with minimal number of DoFs of the form
\begin{align}
	\label{eq:generalVred}
	\V_r^{r-1}(E_N) &= \Curl\,\cDS_{r+1}(E_N)\oplus\x\Po_{r-1}\\\nonumber
	&= \Curl\,\Po_{r+1}(E_N)\oplus\x\Po_{r-1}\oplus\Curl\,\Supp_{r+1}^\cDS(E_N)\\\nonumber
        &= \Po_r^2(E_N)\oplus\Supp_r^\V(E_N),\\
	\label{eq:generalV}
	\V_r^{r}(E_N) &= \Curl\,\cDS_{r+1}(E_N)\oplus\x\Po_{r}\\\nonumber
        &= \Curl\,\Po_{r+1}(E_N)\oplus\x\Po_{r}\oplus\Curl\,\Supp_{r+1}^\cDS(E_N)\\\nonumber
        &= \Po_r^2(E_N)\oplus\x\tilde\Po_{r}\oplus\Supp_r^\V(E_N),
\end{align}
with the following definition of the supplemental (vector valued) functions
\begin{equation}\label{eq:defSuppV}
	\Supp_{r}^{\V}(E_N) = \Curl\,\Supp_{r+1}^{\cDS}(E_N).
\end{equation}

Similar to \cite{Arbogast_Correa_2016,Arbogast_Tao_Wang_2020x_serendipityMixed}, the DoFs ($\cN$ in
Definition~\ref{defn:ciarlet}) for $\v\in\V_r^{s}(E_N)$, $s=r-1,r$, are given (after fixing a basis
for the test functions) by
\begin{alignat}2
\label{eq:vdofEdge}
&\int_{e_i}\v\cdot\nu_i\,p\,d\sigma,&&\quad\forall p\in\Po_r(e_i),\ i=1,2,\ldots,N,\\
\label{eq:vdofDiv}
&\int_{E_N}\v\cdot\grad q\,dx,&&\quad\forall q\in\Po_{s}(E_N),\ q\text{ not constant},\\
\label{eq:vdofCurl}
&\int_{E_N}\v\cdot\bfpsi\,dx,&&\quad\forall \bfpsi\in\Bu_r^\V(E_N),\text{ if }r\ge N-1,
\end{alignat}
where $d\sigma$ is the one dimensional surface measure and the $H^1(E_N)$ and $H(\Div;E_N)$ bubble
functions, for $r\ge N-1$, are
\begin{equation}
\label{eq:bubbles}
\Bu_{r+1}(E_N) = \lambda_1\lambda_2\ldots\lambda_N\Po_{r-N+1}(E_N)
\quad\text{and}\quad
\Bu_r^\V(E_N) = \Curl\,\Bu_{r+1}(E_N).
\end{equation}
We remark that the edge DoFs \eqref{eq:vdofEdge} determine the normal components (flux) of our
vector functions, the divergence DoFs \eqref{eq:vdofDiv} determine the divergence of our vector
functions (with the previous edge DoFs), and the curl DoFs \eqref{eq:vdofCurl} control the curl of
our vector functions.

\begin{theorem}
  The finite element $\V_r^s(E_N)$ defined by \eqref{eq:generalVred}--\eqref{eq:generalV},
  \eqref{eq:defSuppV} for $r=1,2\ldots$ and $s=r-1,r$ (but $s\ge0$) with DoFs
  \eqref{eq:vdofEdge}--\eqref{eq:vdofCurl}, \eqref{eq:bubbles} is well defined (i.e.,
  unisolvent). Moreover, it has the minimal number of DoFs needed of a space of index~$r$ that is $H(\Div)$
  conforming and has independent divergence approximation to order~$s$.
\end{theorem}

\begin{proof}
  The minimal number of DoFs needed are expressed by \eqref{eq:vdofEdge}--\eqref{eq:vdofCurl}, since
  \eqref{eq:vdofEdge} is required for $H(\Div)$ conformity of order~$r$ and \eqref{eq:vdofDiv} is required for
  independent divergence approximation to order~$s$. Moreover, $\eqref{eq:vdofCurl}$ is required to
  control polynomials of degree~$r$ which have no divergence nor edge normal flux.

  The total number of degrees of freedom is
\begin{equation}\label{eq:nDoFsMixed}
  D_{N,r}^\V = \begin{cases}
    N\dim\Po_r(e) + (\dim\Po_s(E_N) - 1) + \dim\Po_{r-N+1}(E_N), &\text{if }r\ge N-1,\\
    N\dim\Po_r(e) + (\dim\Po_s(E_N) - 1), &\text{if }r<N-1,
  \end{cases}
\end{equation}
and the local dimensions of the spaces are
\begin{equation}\label{eq:dimVMixed}
  \dim\V_r^{s}(E_N) = (\dim\cDS_{r+1}-1) + \dim(\x\Po_s).
\end{equation}
  By \eqref{eq:amtotalNumber} and \eqref{eq:totalNumber_small_r}, these numbers agree. In fact,
\begin{equation}\label{eq:dimMixed}
  D_{N,r}^\V = \begin{cases}
  N(r+1) - 1 + \frac12(s+2)(s+1)\\
  \qquad+ \frac12(r-N+3)(r-N+2), & r\ge N-1,\\
  N(r+1) - 1 + \frac12(s+2)(s+1), & r< N-1.
\end{cases}
\end{equation}
The remainder of the proof, to show that these spaces are unisolvent (i.e., a vector function in
$\V^s_r{(E_N)}$ with vanishing DoFs is zero everywhere), is essentially the same as that given in
\cite{Arbogast_Tao_Wang_2020x_serendipityMixed} for direct mixed spaces on quadrilaterals.
\end{proof}


\subsection{Implementation of the mixed method}\label{sec:implementMixed}
The mixed space of vector functions $\V_r^s$ over $\Omega$ is defined by merging continuously the
normal fluxes across each edge~$e$ of the mesh $\cT_h$. That is, for $r\ge0$, $s=r-1,r$, $s\ge0$,
\begin{equation}\label{eq:VinHdiv}
	\V_r^{s} = \big\{\v\in H(\Div;\Omega)\;:\;\v\big|_{E_N}\in\V_r^{s}(E_N)\text{ for all }E_N\in\cT_h\big\}.
\end{equation}
Associated to this space is the scalar space of its divergences, namely,
\begin{equation}
W_s = \Div\,\V_r^s = \big\{w\in L^2(\Omega)\;:\;w\big|_{E_N}\in\Po_s(E_N)\text{ for all }E_N\in\cT_h\big\}.
\end{equation}
It is used, for example, when solving a second order elliptic partial differential equation in mixed form.


\subsubsection{Implementation using the hybrid mixed method}\label{sec:hybrid}

The hybrid form of the mixed method is often used \cite{Arnold_Brezzi_1985} so that no globally
merged basis is required.  A Lagrange multiplier space is used to enforce the normal flux continuity
through an additional equation, using the space
\begin{equation}
\Lambda_r = \big\{\lambda\in L^2\big(\cup_{E_N\in\cT_h}\d E_N\big)\;:\;
                          \lambda\big|_e\in\Po_r(e)\text{ for each edge }e\text{ of }\cT_h\big\}.
\end{equation}

The vector functions in $\V_r^s(E_N)$ can be represented by any of the equivalent forms in
\eqref{eq:generalVred}--\eqref{eq:generalV}.  First, since
$\V_r^s(E_N)=\Curl\,\cDS_{r+1}(E_N)\oplus\x\Po_{s}$, we can construct the full space
$\cDS_{r+1}(E_N)$ as discussed in Sections~\ref{sec:DShigh} and \ref{sec:universal_small_r}, apply the
curl operator, and add in $\x\Po_s(E_N)$. But we can also use the fact that
$\V_r^{r-1}(E_N)=\Po_r^2(E_N)\oplus\Supp_r^\V(E_N)$ and
$\V_r^{r}(E_N)=\Po_r^2(E_N)\oplus\x\tilde\Po_{r}\oplus\Supp_r^\V(E_N)$, and simply add to the
polynomials the supplemental space $\Supp_{r}^{\V}(E_N)=\Curl\,\Supp_{r+1}^{\cDS}(E_N)$. To
construct $\Supp_{r+1}^{\cDS}(E_N)$, one uses \eqref{eq:supplementFcns}--\eqref{eq:supplementSpace}
when $r$ is large, and otherwise requires the construction given in
Section~\ref{sec:whole_space_small_r}.


\subsubsection{Implementation as an $H(\Div)$-conforming mixed space}\label{sec:conforming}

If an explicit basis for the $H(\Div)$-conforming space \eqref{eq:VinHdiv} of vector-valued
functions is required, one can proceed as follows. The construction is an extension of the $N=4$
case given in \cite{Arbogast_Tao_Wang_2020x_serendipityMixed}. We use the fact that the tangential
derivative of a function along an edge $e_i$ of an element $E_N$ maps by the $\Curl$ operator to a
normal derivative, i.e., for $\phi\in\cDS_{r+1}(E_N)$,
\begin{equation}\label{eq:tangentialToNormalDer}
  \grad\phi\cdot\tau_i\big|_{e_i}
  = \Curl\,\phi\cdot\nu_i\big|_{e_i}, \quad\text{with } \tau_i=(-\nu_{i,2},\nu_{i,1})\text { on }e_i.
\end{equation}
Since the serendipity spaces are globally continuous, the tangential derivatives will agree
across~$e_i$, which implies that the global basis functions arising from $\cDS_{r+1}(\Omega)$ will
be in $H(\Div;\Omega)$.

We construct $H(\Div)$-conforming vector basis functions in four sets, related to the edge DoFs
\eqref{eq:vdofEdge} with nonconstant test functions, the edge DoFs \eqref{eq:vdofEdge} with constant
test functions, the divergence DoFs \eqref{eq:vdofDiv}, and the curl DoFs \eqref{eq:vdofCurl}.

\paragraph{Basis functions from curls of interior cell basis functions of $\cDS_{r+1}(E_N)$}
The interior cell basis functions of $\cDS_{r+1}(E_N)$ are
$\{\varphi_{E,i}^{(r+1)},\ i=1,2,\ldots,\dim\Po_{r+1-N}\}$ as given by \eqref{eq:nodal-cell} (the
superscript is a reminder that the index of the direct serendipity space is $r+1$). However, any
basis for \eqref{eq:nodal-cell-space}, i.e., the bubble space $\Bu_{r+1}(E_N)$ defined in
\eqref{eq:bubbles}, suffices.  Denote it as
$\{\phi_{E_N,i}^{(r+1)},\ i=1,2,\ldots,\dim\Po_{r+1-N}\}$.  Then for each $E_N\in\cT_h$, the global
basis functions for $\V_r^s$ are
\begin{equation}\label{eq:mixed-basis-bubble}
  \bfpsi_{b,E_N,i} = \begin{cases}\Curl\,\phi_{E_N,i}^{(r+1)},\quad i=1,\ldots,\tfrac12(r+3-N)(r+2-N),
    &\text{on }E_N,\\
    0,&\text{otherwise}.
  \end{cases}
\end{equation}
These exist only when $r\ge N-1$, and they are in fact the $H(\Div)$ bubble functions $\Bu_r^{\V}$
appearing in \eqref{eq:bubbles}.  They have no normal flux and no divergence. They are associated to
the curl DoFs \eqref{eq:vdofCurl}.

\paragraph{Basis functions from curls of interior edge basis functions of $\cDS_{r+1}(E_N)$}
The interior edge basis functions of $\cDS_{r+1}(E_N)$ are
$\{\varphi_{e,i,j}^{(r+1)},\ i=1,2,\ldots,N,\ j=1,2,\ldots,r\}$ as given by \eqref{eq:nodal-edge11}
or \eqref{eq:varphi_r_e11} when $r<N-2$.  For $r\ge N-2$, one could use the simpler set
$\{\phi_{e,i,j}^{(r+1)}/\phi_{e,i,j}^{(r+1)}(\x_{e,i,j})\}$ given in \eqref{eq:nodal-edge-tilde}
which ignores the internal cell DoFs, and we proceed with this choice (the case $r<N-2$ is entirely
similar).  Consider an edge $e$ of the mesh shared by elements $E_{k}$ and $E_{\ell}$ with $k<\ell$
and $e$ locally denoted as edge $i_1$ and $i_2$, respectively. The global basis functions for
$\V_r^s$ are, for $r\ge 1$ and $j=1,\ldots,r$,
\begin{equation}\label{eq:mixed-basis-interiorEdge}
  \bfpsi_{e,j}(\x) = \begin{cases}
    \Curl\,\phi_{e,i_1,j}^{(r+1)}(\x)/\phi_{e,i_1,j}^{(r+1)}(\x_{e,i_1,j}),&\x\in E_{k},\\[6pt]
    \Curl\,\phi_{e,i_2,r-j+1}^{(r+1)}(\x)/\phi_{e,i_2,r-j+1}^{(r+1)}(\x_{e,i_2,r-j+1}),&\x\in E_\ell,\\
    0,&\x\notin E_{k}\cup E_\ell.
\end{cases}
\end{equation}
These functions have vanishing divergence but nonvanishing normal flux; however, the average normal
flux vanishes. They are associated to the edge DoFs \eqref{eq:vdofEdge} with nonconstant test functions.

\paragraph{Basis functions from curls of vertex basis functions of $\cDS_{r+1}(E_N)$}
We will now construct basis functions that have constant normal flux on a single edge of the
mesh. These cannot have vanishing divergence. We will use the vertex basis functions of
$\cDS_{r+1}(E_N)$, which are $\{\varphi_{v,i}^{(r+1)},\ i=1,2,\ldots,N\}$ as given in
\eqref{eq:nodal-vertices} or \eqref{eq:varphi_r_v1}.  Again, when $r\ge N-2$ we can instead simply
use $\{\phi_{v,i}^{(r+1)}/\phi_{v,i}^{(r+1)}(\x_{v,i})\}$ given in
\eqref{eq:nodal-vertices-unnormalized}, and we proceed with the discussion using this case. The
construction is complicated by the fact that the curls of these functions have nonvanishing normal
flux on all the edges of the mesh emanating from the vertex in question.

We work on the element $E_N$, and we first modify the serendipity vertex basis functions so that their
restrictions to each edge $e$ of $E_N$ is a linear function, i.e., we define for all~$i$
\begin{align*}
	\phi^*_{v,i}(\x) = \frac{\phi^{(r+1)}_{v,i}(\x)}{\phi^{(r+1)}_{v,i}(\x_{v,i})} + \sum_{j=1}^r \bigg[
	\frac{j}{r+1}\frac{\phi^{(r+1)}_{e,i,j}(\x)}{\phi^{(r+1)}_{e,i,j}(\x_{e,i,j})} + \Big(1 - \frac{j}{r+1}\Big) \frac{\phi^{(r+1)}_{e,i+1,j}(\x)}{\phi^{(r+1)}_{e,i+1,j}(\x _{e,i+1,j})}\bigg],
\end{align*}
again using indices modulo~$N$. Then define $\bfpsi^*_{v,i}=\Curl\,\phi^*_{v,i}$, for which
\begin{equation*}
  \bfpsi^*_{v,i}(\x)\cdot\nu_j\big|_{e_j} = \grad\phi^*_{v,i}(\x)\cdot\tau_j\big|_{e_j}
  = \begin{cases}
    1/|e_i|,&j=i,\\
    -1/|e_{i+1}|,&j=i+1,\\
    0,&\text{otherwise}.
  \end{cases}
\end{equation*}
We also use the vector $\bfpsi^{**}_{v,i}(\x) = \x-\x_{v,i+1}\in\x\Po_0(E_N)\oplus\Po_0^2(E_N)\subset\V_r^s(E_N)$, which is in
our space and satisfies
\begin{equation*}
  \bfpsi^{**}_{v,i}(\x)\cdot\nu_j\big|_{e_j} 
   = \begin{cases}
    0,&j=i+1,i+2, \\
    (\x_{v,j}-\x_{v,i+1})\cdot\nu_j,&\text{otherwise},
    
  \end{cases}
\end{equation*}
which is nonnegative on every edge $e_j$.

For any edge $e_i$ of element $E_N$, we define a vector function with flux only on $e_i$ by
canceling the fluxes of $\bfpsi^{**}_{v,i}$ on all the other edges using some of the
$\bfpsi^*_{v,k}$. Precisely, we define for edge $e=e_i$ of element $E_N$
\begin{gather}\label{eq:mixed-basis-constantEdge-STAR}
  \bfpsi_{e,0}\big|_{E_N}
  = \frac{1}{c_{i,i+N}}\Bigg(\bfpsi^{**}_{v,i}-\sum_{j=i+3}^{i+N-1}c_{i,j}\,|e_{j}|\,\bfpsi^*_{v,j}\Bigg), \\\nonumber
  c_{i,i+2} = 0,\quad c_{i,j} = (\x_{v,j}-\x_{v,i+1})\cdot\nu_{j} + \frac{|e_{j-1}|}{|e_{j}|}\,c_{i,j-1} > 0,\quad
  j = i+3,\ldots,i+N,
\end{gather}
which has normal flux 1 on $e_i$ and 0 on all the other edges. These can be merged across edges to
define $H(\Div)$-conforming global basis functions, which have constant divergence on each
element. Note that the choice of vertex index $\x_{v,i+1}$ in $\bfpsi^{**}_{v,i}$ is only for
convenience in presenting the construction. We might have chosen it to be any other vertex except
$\x_{v,i-1}$ and $\x_{v,i}$. The basis functions here are associated to the edge DoFs
\eqref{eq:vdofEdge} with constant test functions.

\paragraph{Basis functions with nonvanishing and nonconstant divergence}
Finally, when $s\ge1$ we define the global basis functions associated to the nonconstant
divergences.  They are local to each element $E_N\in\cT_h$. Working on $E_N$, we begin with the
functions $\x\Po_s^*(E_N)$, where $\Po_s^*(E_N)=\sum_{k=1}^s\tilde\Po_k(E_N)\subset\Po_s(E_N)$. Take
$p_i(\x)$ in a basis for $\Po_s^*(E_N)$, so $i=1,\ldots,\tfrac12(s+2)(s+1)-1$.  We must remove the
normal flux on $\d E_N$ from $\x p_i(\x)$.  We do this using \eqref{eq:mixed-basis-interiorEdge} and
\eqref{eq:mixed-basis-constantEdge} by defining
\begin{equation}\label{eq:mixed-basis-divergence}
  \bfpsi_{d,E_N,i}(\x) = \begin{cases}\displaystyle
    \x p_i(\x) - \sum_{j=1}^N\sum_{k=0}^r\alpha_{j,k}\,\bfpsi_{e_j,k}(\x),&\text{on }E_N,\\
    0,&\text{otherwise},
    \end{cases}
\end{equation}
and setting the coefficients $\alpha_{j,k}$ on each edge $e_j$ so that
\begin{equation}\label{eq:mixed-basis-divergence-alpha}
0 =c_j p_i(\x) - \sum_{k=0}^r\alpha_{j,k}\,\bfpsi_{e_j,k}(\x)\cdot\nu_j\big|_{e_j},
\end{equation}
where $c_j=\x\cdot\nu_j|_{e_j}$ is a constant.  The coefficients can be found once one realizes that
on edge $e_j$, $\varphi_{e,j,k}^{(r+1)}(x)\cdot\tau_j\big|_{e_j}=\mathfrak{L}_k(t)$, a Lagrange
basis polynomial, where $\x(t)=(1-t)\x_{v,j-1}+t\x_{v,j}$ for $t\in[0,1]$.  Therefore, for $k\ge1$,
\begin{equation*}
  \bfpsi_{e_j,k}(\x)\cdot\nu_j\big|_{e_j}
  = \Curl\,\varphi_{e,j,k}^{(r+1)}(x)\cdot\nu_j\big|_{e_j}
  = \grad\varphi_{e,j,k}^{(r+1)}(x)\cdot\tau_j\big|_{e_j}
  = \frac{\mathfrak{L}_k'(t)}{|\x_{v,j} - \x_{v,j-1}|},
\end{equation*}
and
\begin{equation}
0 = c_j \int_0^tp_i(\x(s))\,ds - \alpha_{j,0}\,t - \sum_{k=1}^r\frac{\alpha_{j,k}\,\mathfrak{L}_k(t)}{|\x_{v,j} - \x_{v,j-1}|}.
\end{equation}
The coefficients can be read off by substituting in the Lagrange points $t_\ell=\ell/(r+1)$ for
$\ell=1,\ldots,r+1$.  These basis functions are associated to the divergence DoFs \eqref{eq:vdofDiv}
with nonconstant local divergence.

The global basis is now fully defined.


\section{Approximation properties for $\V_r^s$}\label{sec:approxMixed}

In this section, we state the approximation theory for our new direct mixed finite elements. A
discussion and detailed proof for the $N=4$ case has been given in
\cite{Arbogast_Tao_Wang_2020x_serendipityMixed}. The proof for polygons is very similar, and so
omitted here.

We can define a projection operator $\pi:H(\Div;\Omega)\cap(L^{2+\epsilon}(\Omega))^2\to\V_r^s$,
$s=r-1,r$, where $\epsilon>0$, by piecing together locally defined operators $\pi_E$. For suitable
$\v$, $\pi_E\v$ is defined in terms of the DoFs \eqref{eq:vdofEdge}--\eqref{eq:vdofCurl}. The
operator $\pi$ satisfies the commuting diagram property~\cite{Douglas_Roberts_1985}, which is to say
that
\begin{equation}
\cP_{W_s}\div\v=\div\pi\v,
\end{equation}
where $\cP_{W_s}$ is the $L^2$-orthogonal projection operator onto $W_s=\div\V_r^s$. The following
lemma holds.

\begin{theorem}\label{thm:brambleMixed}
  Let $\cT_h$ be uniformly shape regular with shape regularity parameter $\sigma_*$ and let
  Assumption~\ref{assumption:niceDS} hold.  Then for $\V_r^s$ there is a constant
  $C=C(r,\sigma^*)>0$, independent of $h>0$, such that
\begin{alignat}2
\label{eq:approx_u}
\|\v-\pi\v\|_{L^2(\Omega)} &\le C\,\|\v\|_{H^k(\Omega)}\,h^{k},&&\quad k=1,\ldots,r+1,\\
\label{eq:approx_p}
\|p-\cP_{W_s}p\|_{L^2(\Omega)} &\le C\,\|p\|_{H^k(\Omega)}\,h^{k},&&\quad k=0,1,\ldots,s+1,\\
\label{eq:approx_divu}
\|\div(\v-\pi\v)\|_{L^2(\Omega)} &\le C\,\|\div\v\|_{H^k(\Omega)}\,h^{k},&&\quad k=0,1,\ldots,s+1,
\end{alignat}
where $s=r-1\ge0$ and $s=r\ge1$ for reduced and full $H(\Div)$-approximation, respectively.
Moreover, the discrete inf-sup condition
\begin{equation}
\label{eq:inf-sup}
\sup_{\v_h\in\V_r^s}\frac{(w_h,\div\v_h)}{\|\v_h\|_{H(\text{\rm div})}}
 \ge \gamma\,\|w_h\|_{L^2(\Omega)},\quad\forall w_h\in W_s,
\end{equation}
holds for some $\gamma=\gamma(r,\sigma^*)>0$ independent of $h>0$.
\end{theorem}


\section{Numerical results}\label{sec:numerics}

We test our finite elements on Poisson's equation
\begin{alignat}{2}
\label{eq:elliptic-pde}
-\div(\grad p) &= f &&\quad\text{in }\Omega,\\
\label{eq:elliptic-dirBC}
 p &= 0 &&\quad\text{on }\partial\Omega,
\end{alignat}
where
$f\in L^2(\Omega)$.  The problem can be written in the weak form: Find
$p\in H_0^1(\Omega)$ such that
\begin{align}
\label{eq:weak-bvp}
(\grad p,\grad q) = (f,q), \quad\forall q\in H_0^1(\Omega),
\end{align}
where $(\cdot,\cdot)$ is the $L^2(\Omega)$ inner product. Setting
\begin{equation}
\label{eq:elliptic-flux}
\u = -\grad p,
\end{equation}
we also have the mixed weak form: Find $\u\in H(\Div;\Omega)$ and $p\in L^2(\Omega)$ such that
\begin{alignat}3
\label{eq:mixed-darcy}
&(\u,\v) - (p,\div\v) &&=0, &&\quad\forall\v\in H(\Div;\Omega),\\
\label{eq:mixed-conservation}
&(\div\u,w) &&= (f,w), &&\quad\forall w\in L^2(\Omega).
\end{alignat}
These weak forms give rise to finite element approximations. In view of Theorems~\ref{thm:approxDS}
and~\ref{thm:brambleMixed}, it is well known that the following theorem holds
\cite{Brezzi_Fortin_1991, Brenner_Scott_1994}.

\begin{theorem}\label{thm:convergence}
  Let $\cT_h$ be uniformly shape regular with shape regularity parameter $\sigma_*$ and let
  Assumption~\ref{assumption:niceDS} hold. There exists a constant $C>0$, depending on $r$ and
  $\sigma_*$ but otherwise independent of $\cT_h$ and $h>0$, such that
\begin{align}
\|p-p_h\|_{H^m(\Omega)}&\leq C\,h^{s+1-m}\,|p|_{H^{s+1}(\Omega)},\quad s=0,1,\ldots,r,\quad m=0,1,
\end{align}
where $p_h\in\cDS_r(\Omega)\cap H_0^1(\Omega)$ approximates \eqref{eq:weak-bvp} for
$r\ge1$. Moreover,
\begin{alignat}2
\|\u-\u_h\|_{L^2(\Omega)} &\le C\|\u\|_{H^k(\Omega)}h^{k},&&\quad k=1,\ldots,r+1,\\
\|p-p_h\|_{L^2(\Omega)} &\le C\|\u\|_{H^k(\Omega)}h^{k},&&\quad k=1,\ldots,s+1,\\
\|\div(\u-\u_h)\|_{L^2(\Omega)} &\le C\|\div\u\|_{H^k(\Omega)}h^{k},&&\quad k=0,1,\ldots,s+1,
\end{alignat}
where $(\u_h,p_h)\in\V_r^s\times W_s$ approximates
\eqref{eq:mixed-darcy}--\eqref{eq:mixed-conservation}, for $r\ge0$ and $0\le s=r,r-1$.
\end{theorem}

We consider the test problem~\eqref{eq:elliptic-pde}--\eqref{eq:elliptic-dirBC}
defined on the unit square $\Omega = [0,1]^2$.  The exact solution is
$u(x_1,x_2) = \sin(\pi x_1)\sin(\pi x_2)$ and the source term is $f(\x) = 2\pi^2\sin(\pi x_1)\sin(\pi x_2)$.

\begin{figure}[ht]
\centerline{
  \parbox{.2\linewidth}{\includegraphics[width=\linewidth]{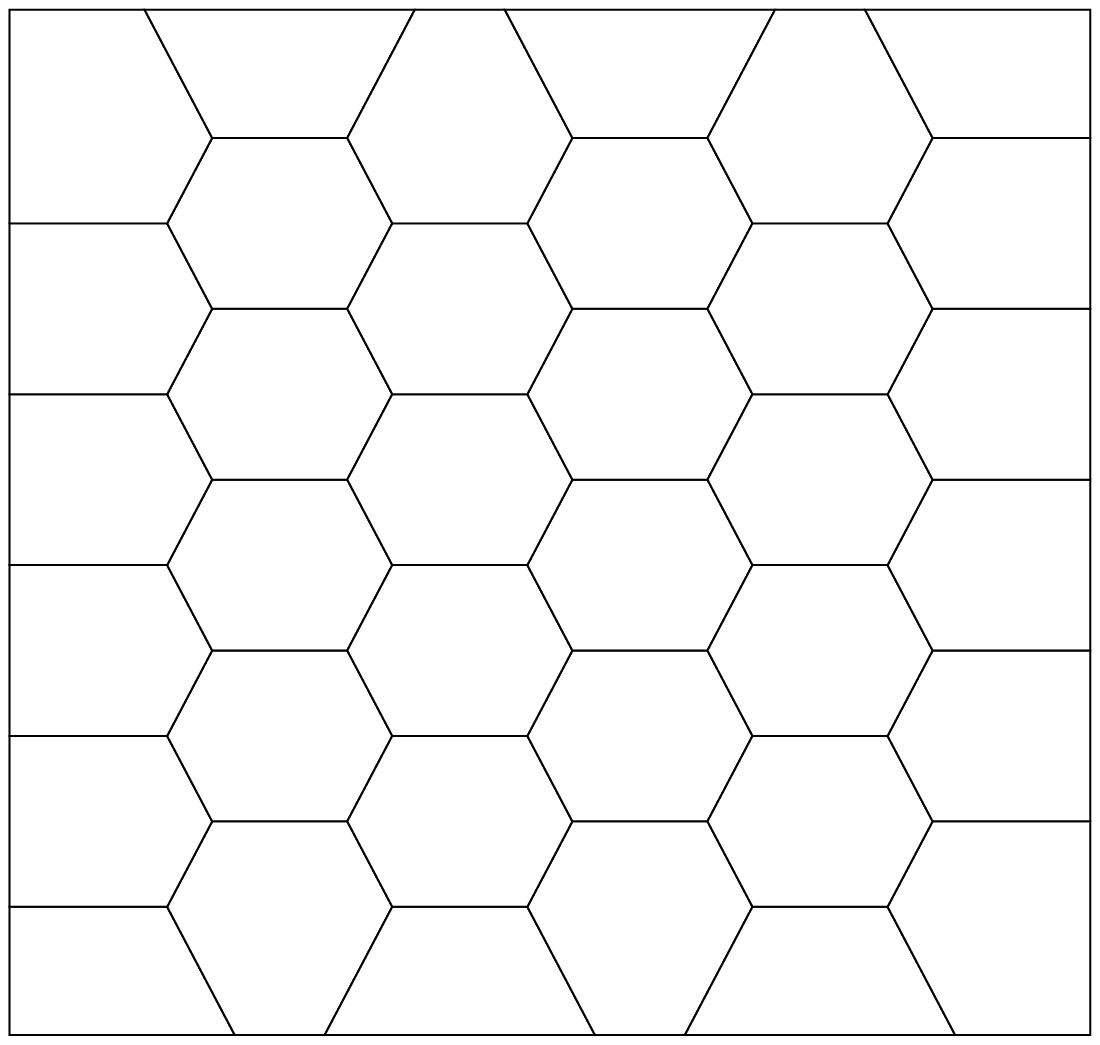}\\[-4pt]
    \centerline{\scriptsize$\cT_h^1$, $n=6$}}
  \quad\parbox{.2\linewidth}{\includegraphics[width=\linewidth]{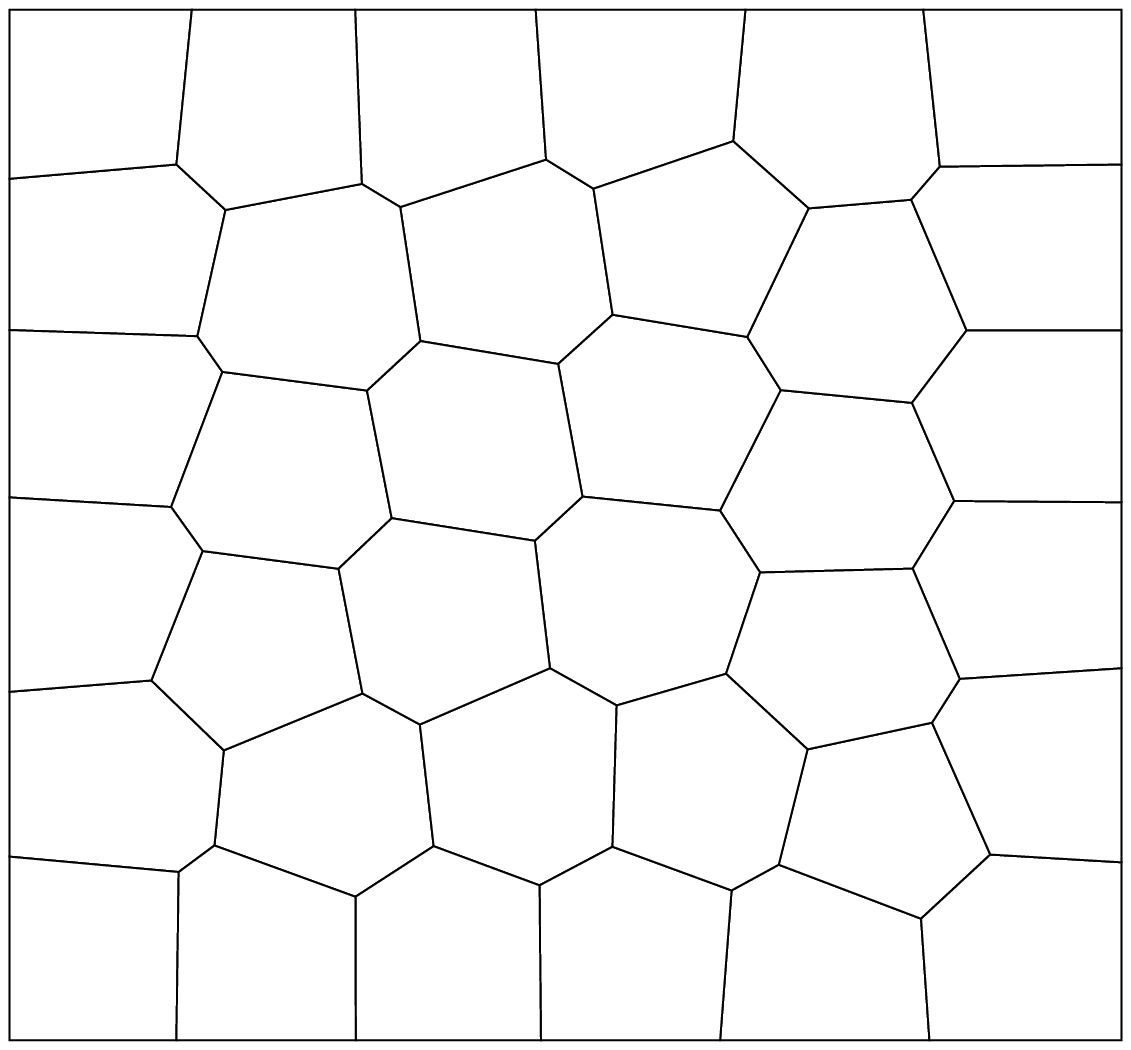}\\[-4pt]
    \centerline{\scriptsize$\cT_h^2$, $n=6$}}
  \quad\parbox{.2\linewidth}{\includegraphics[width=\linewidth]{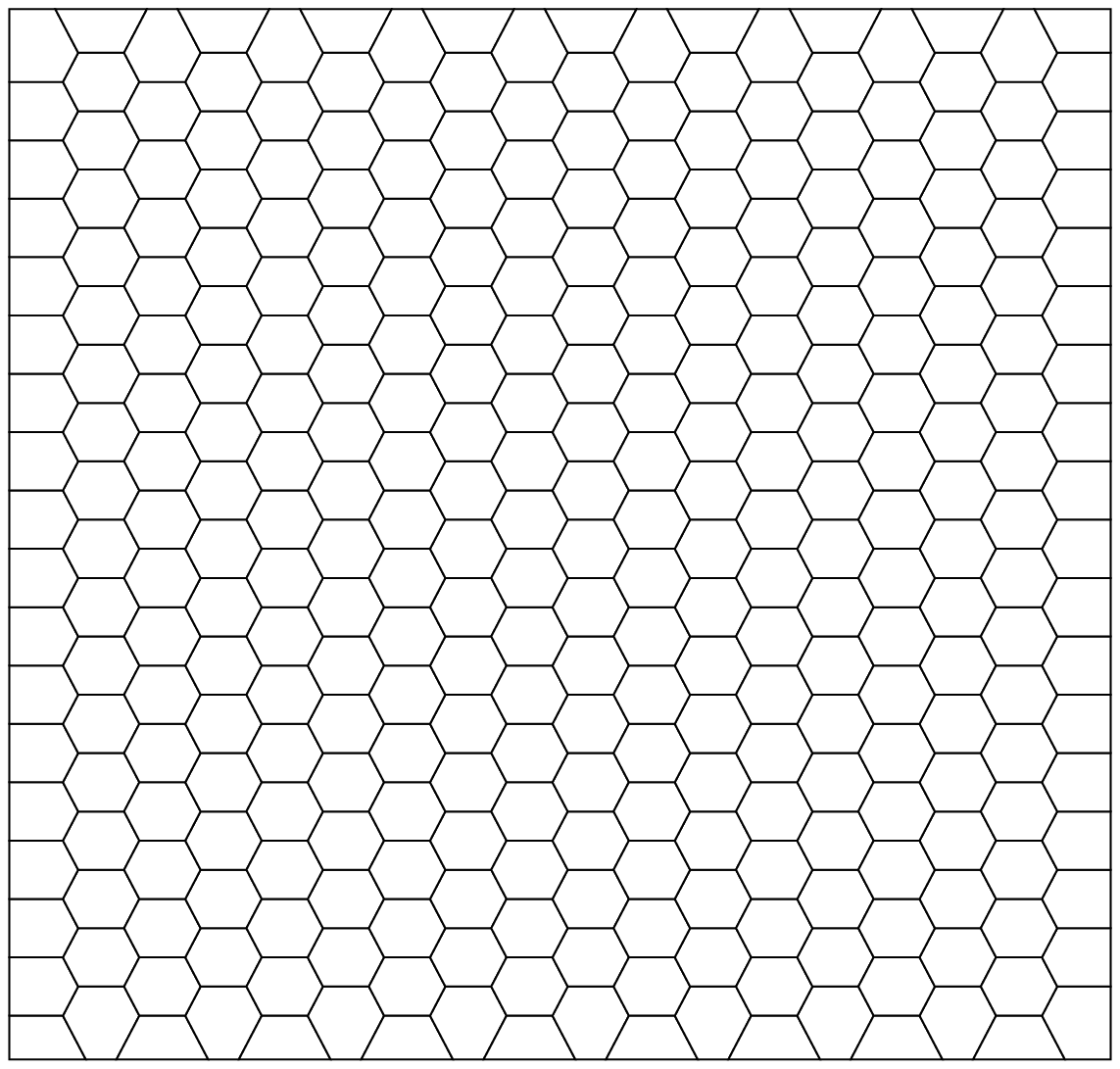}\\[-4pt]
    \centerline{\scriptsize$\cT_h^1$, $n=18$}}
  \quad\parbox{.2\linewidth}{\includegraphics[width=\linewidth]{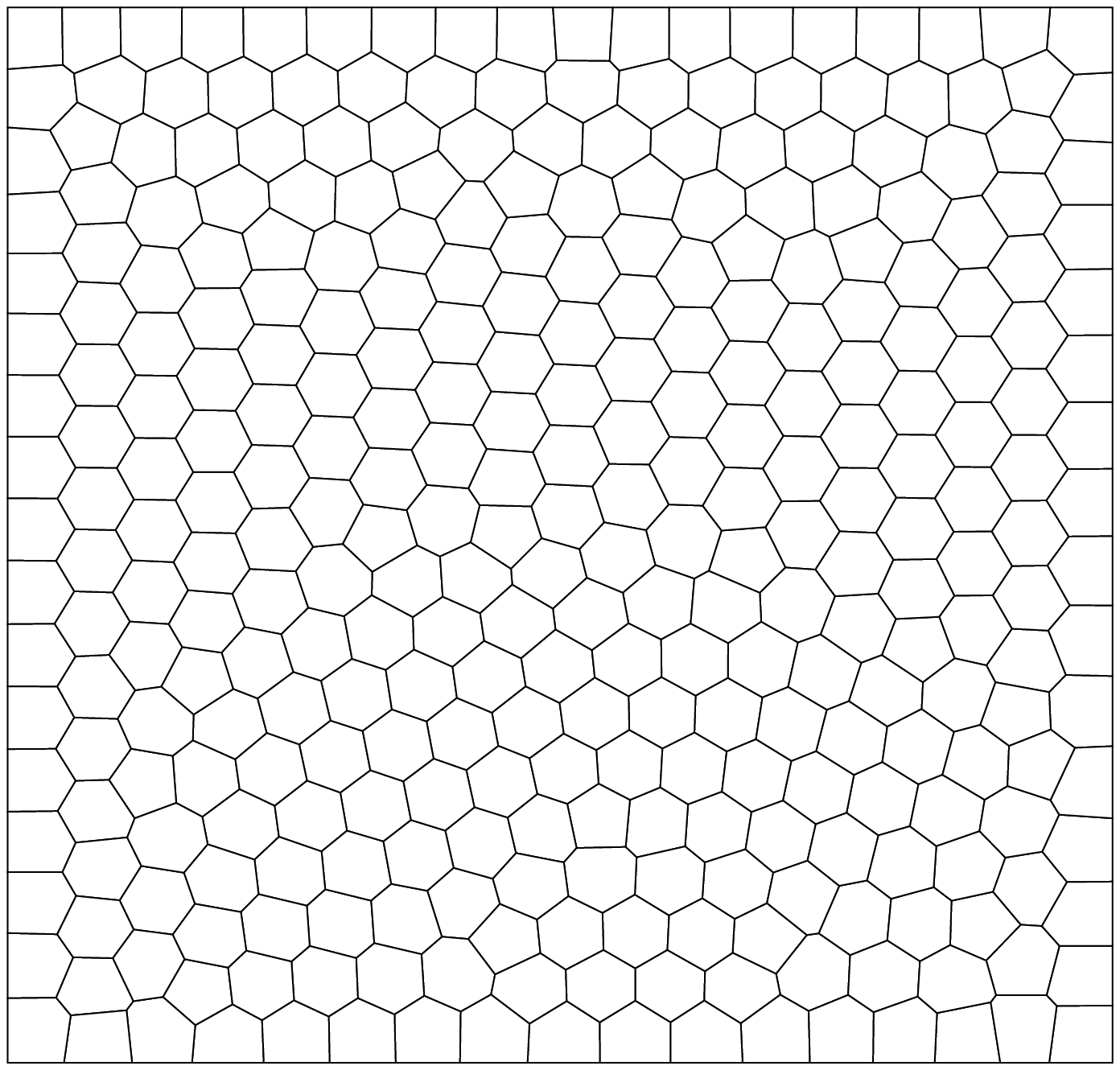}\\[-4pt]
    \centerline{\scriptsize $\cT_h^2$, $n=18$}}
}
\caption{Meshes with $6\times6$ and $18\times18$ elements.}\label{fig:meshes}
\end{figure}

Solutions are computed on two different sequences of meshes, each has $n^2$ elements and is a
Voronoi diagram mesh generated using the software package PolyMesher \cite{TPPM_2012_polymesher}.
The first set of meshes, $\cT_h^1$, is a simple mesh composed of polygons generated from regularly
spaced seeds.  The seeds are initially uniformly spaced and then alternatively perturbed up or down
in the $y$-direction by one quarter of the regular spacing. The number of vertices of each element
is $N=4$, $5$, or $6$. The second sequence, $\cT_h^2$, is generated by PolyMesher using $n^2$ random
initial seeds and up to $10,000$ iterations to smooth the mesh. We illustrate these patterns by
showing the $n = 6$ and $n = 18$ cases in Figure~\ref{fig:meshes}.

\begin{table}[ht]
\caption{Maximum, minimum, and average shape regularity parameters for each mesh.}\label{tab:chunk}
\vspace*{-4pt}\begin{center}\small
\begin{tabular}{c|ccc|ccc|ccc}
\hline 
& \multicolumn{3}{c|}{$\cT_h^1$} 
  & \multicolumn{3}{c|}{$\cT_h^2$} &  \multicolumn{3}{c}{Modified $\cT_h^2$ ($n=18,22$)}\\
$n$ & max & min & average & max & min & average & max & min & average\\
\hline
\06 & 0.568 & 0.355 & 0.401 & 0.778 & 0.180 & 0.341 & --\!--- & --\!--- & --\!---\\
10  & 0.568 & 0.355 & 0.391 & 0.762 & 0.115 & 0.381 & --\!--- & --\!--- & --\!---\\
14  & 0.568 & 0.355 & 0.387 & 0.787 & 0.161 & 0.408 & --\!--- & --\!--- & --\!---\\
18  & 0.568 & 0.355 & 0.384 & 0.787 & 0.127 & 0.378 & 0.787 & 0.160 & 0.380 \\
22  & 0.568 & 0.355 & 0.383 & 0.783 & 0.150 & 0.386 & 0.776 & 0.186 & 0.390 \\
\hline
\end{tabular}
\end{center}
\end{table}

We give results on each mesh sequence for $n=6$, 10, 14, 18, and 22. The maximum, minimum, and average
shape regularity parameters are shown in Table~\ref{tab:chunk}. Sequence $\cT_h^1$ has a
fixed maximum and minimum shape regularity parameter; moreover, the average shape regularity
parameter decreases and converges to a constant as the number of elements increases. However, since
the meshes of $\cT_h^2$ are generated randomly, we can see in Figure~\ref{fig:meshes} that there is
no fixed pattern in the shape of the elements, and so the shape regularity parameter varies as well.
The $n = 18$ and $22$ meshes seem to be less regular than the other $\cT_h^2$ meshes, so to improve
the regularity, we removed some of the small edges, creating the ``modified $\cT_h^2$'' mesh
sequence, as described later in Section~\ref{sec:DS_nonregularMesh}.


\subsection{Direct serendipity spaces}

We present in this section convergence studies for the direct serendipity spaces $\cDS_r$.


\subsubsection{Shape regular meshes of mostly hexagons, $\cT_h^1$}

Table~\ref{tab:t2} shows the errors and orders of convergence for the mesh sequence $\cT_h^1$
consisting of quadrilaterals, pentagons, and hexagons. The convergence rates are consistent with the
theory.

We observed (in results not reported here) that for the same number of elements, the error on a mesh
from $\cT_h^1$ is smaller compared to a mesh of trapezoids. As $n$ increases, the $\cT_h^1$ meshes
are refined in a fixed pattern, giving a higher percentage of elements that are hexagons in the
interior of the mesh. This observation suggests that elements with more edges might tend to give
better approximations.

\begin{table}[ht]
\caption{Errors and convergence rates for $\cDS_r$ on $\cT^1_h$ meshes.}\label{tab:t2}
\vspace*{-4pt}\begin{center}\small
\begin{tabular}{c|cc|cc|cc|cc}
\hline & \multicolumn{2}{c|}{$r=2$} & \multicolumn{2}{c|}{$r=3$}
& \multicolumn{2}{c|}{$r=4$} & \multicolumn{2}{c}{$r=5$} \\
$n$ & error & rate & error & rate & error & rate & error & rate \\
\hline
\multicolumn{9}{c}{$L^2$-errors and convergence rates}\\
\hline
10 & 1.991e-04 & 3.19 & 8.639e-06 & 4.31 & 3.549e-07 & 5.37 & 9.891e-09 & 6.50 \\
14 & 6.960e-05 & 3.12 & 2.129e-06 & 4.16 & 5.921e-08 & 5.32 & 1.152e-09 & 6.39 \\
18 & 3.199e-05 & 3.09 & 7.595e-07 & 4.10 & 1.568e-08 & 5.29 & 2.384e-10 & 6.27 \\
22 & 1.725e-05 & 3.08 &	3.357e-07 &	4.07 &	5.460e-09 &	5.26 &	6.442e-11 &	6.52 \\

\hline
\multicolumn{9}{c}{$H^1$-seminorm errors and convergence rates}\\
\hline
10 & 3.223e-03 & 2.18 & 1.826e-04 & 3.19 & 8.844e-06 & 4.34 & 2.669e-07 & 5.44 \\
14 & 1.575e-03 & 2.13 & 6.441e-05 & 3.10 & 2.083e-06 & 4.30 & 4.383e-08 & 5.37 \\
18 & 9.285e-04 & 2.10 & 2.985e-05 & 3.06 & 7.138e-07 & 4.26 & 1.150e-08 & 5.32 \\
22 & 6.110e-04 & 2.09 &	1.622e-05 &	3.04 &	3.052e-07 &	4.23 &	3.978e-09 &	5.29 \\
\hline
\end{tabular}
\end{center}
\end{table}

\begin{figure}[ht]
\centerline{
\parbox{.5\linewidth}{\centerline{\includegraphics[width=0.60\linewidth]{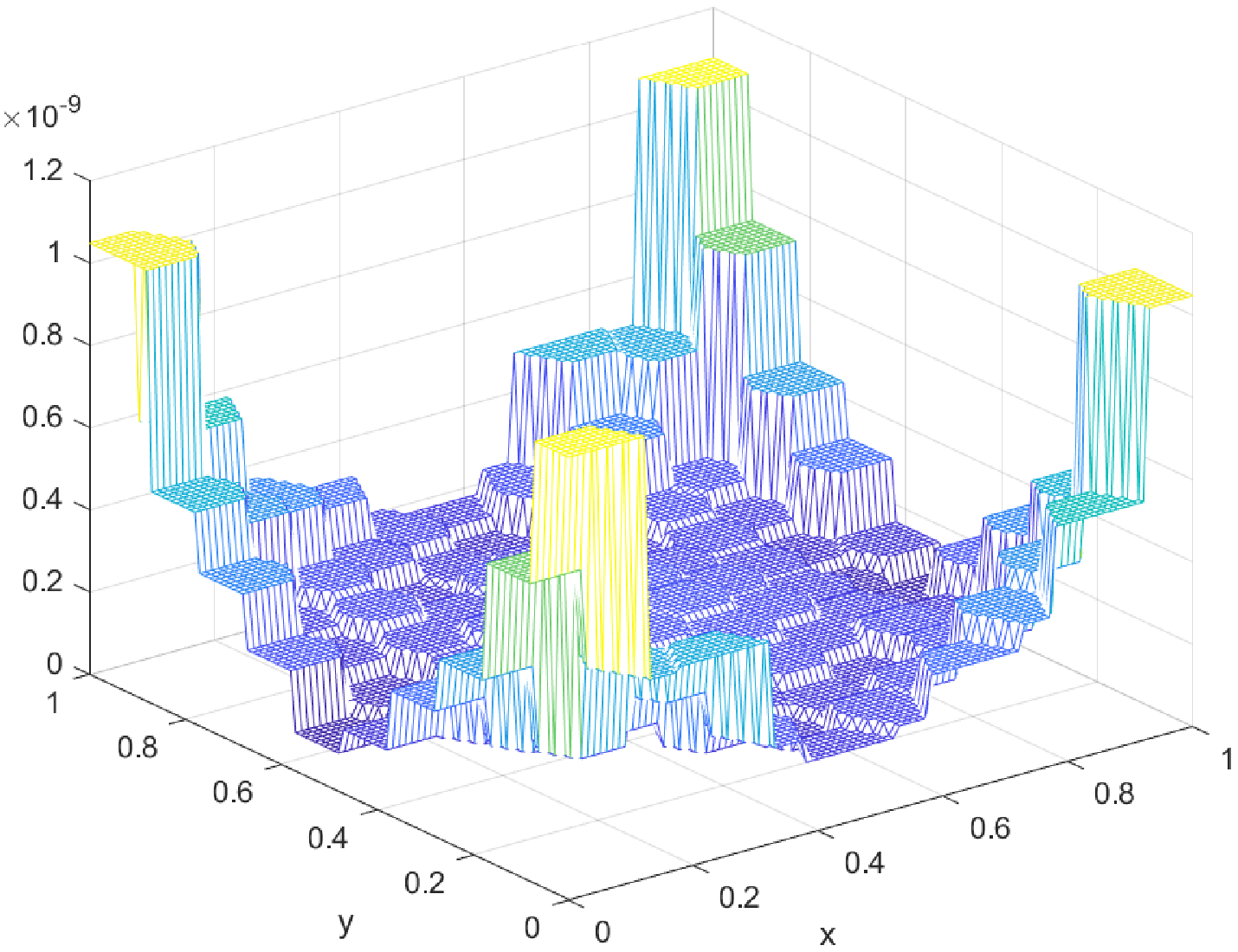}}\centerline{\scriptsize$n=10$}}
\parbox{.5\linewidth}{\centerline{\includegraphics[width=0.60\linewidth]{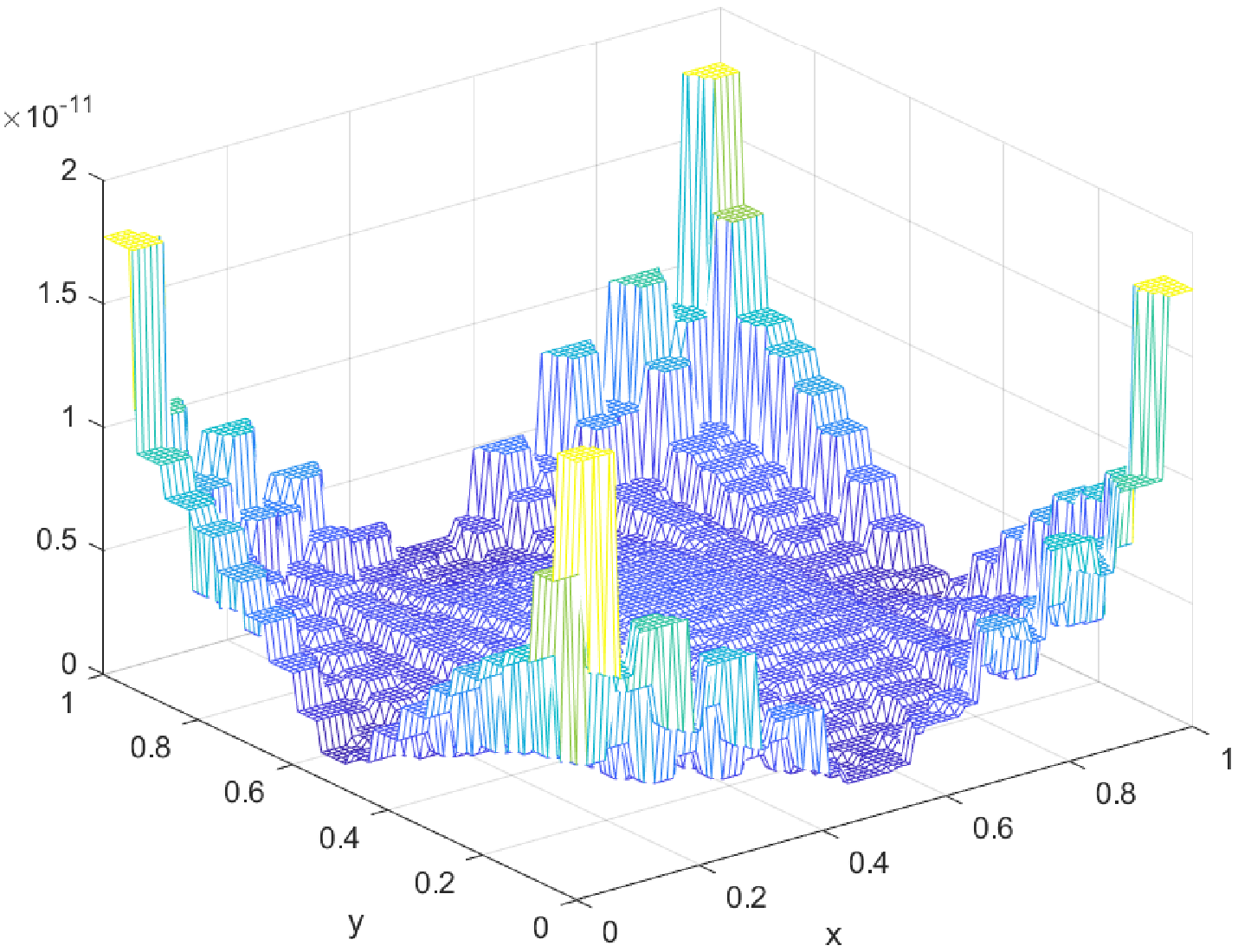}}\centerline{\scriptsize$n=18$}}
}
\caption{The $L^2$ error on each element for mesh sequence $\cT^1_h$ at level $n = 10$ and $n = 18$ with
  approximation index $r=5$.}\label{fig:error_t2}
\end{figure}
	
To test this hypothesis, we graphed the $L^2$ error on each element in Figure~\ref{fig:error_t2} at
level $n = 10$ and 18 with $r = 5$.  The error is indeed concentrated around the boundary, where the
quadrilateral and pentagonal elements concentrate. However, the solution
$u(x_1,x_2) = \sin(\pi x_1)\sin(\pi x_2)$ on $[0,1]^2$ has a single hump over the domain, so the
solution is steepest near the boundary and thus harder to approximate there.

\begin{figure}[ht]
\centerline{
\parbox{.5\linewidth}{\centerline{\includegraphics[width=0.60\linewidth]{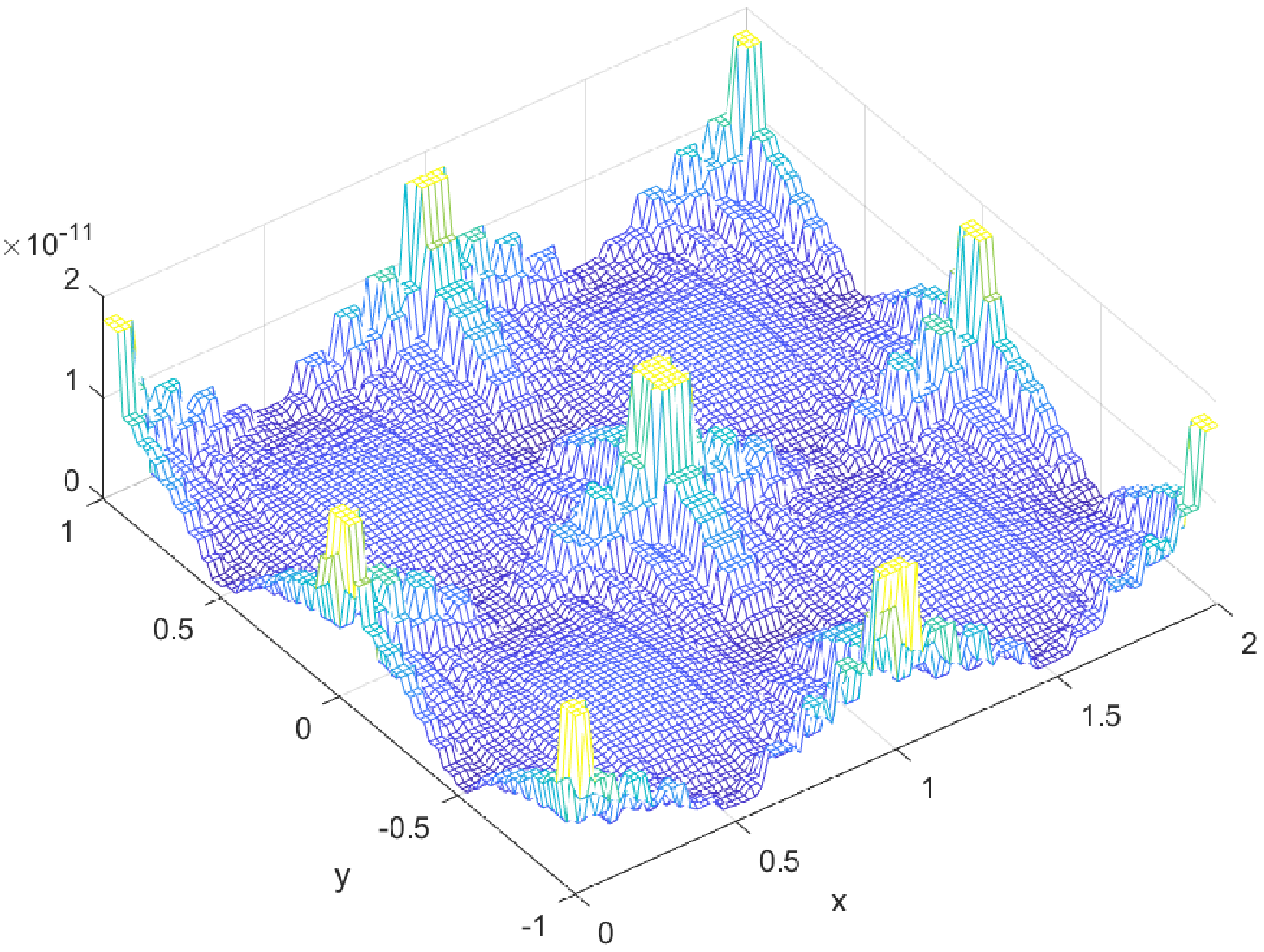}}\vspace*{-4pt}\centerline{\scriptsize Larger domain}}
\parbox{.5\linewidth}{\centerline{\includegraphics[width=0.60\linewidth]{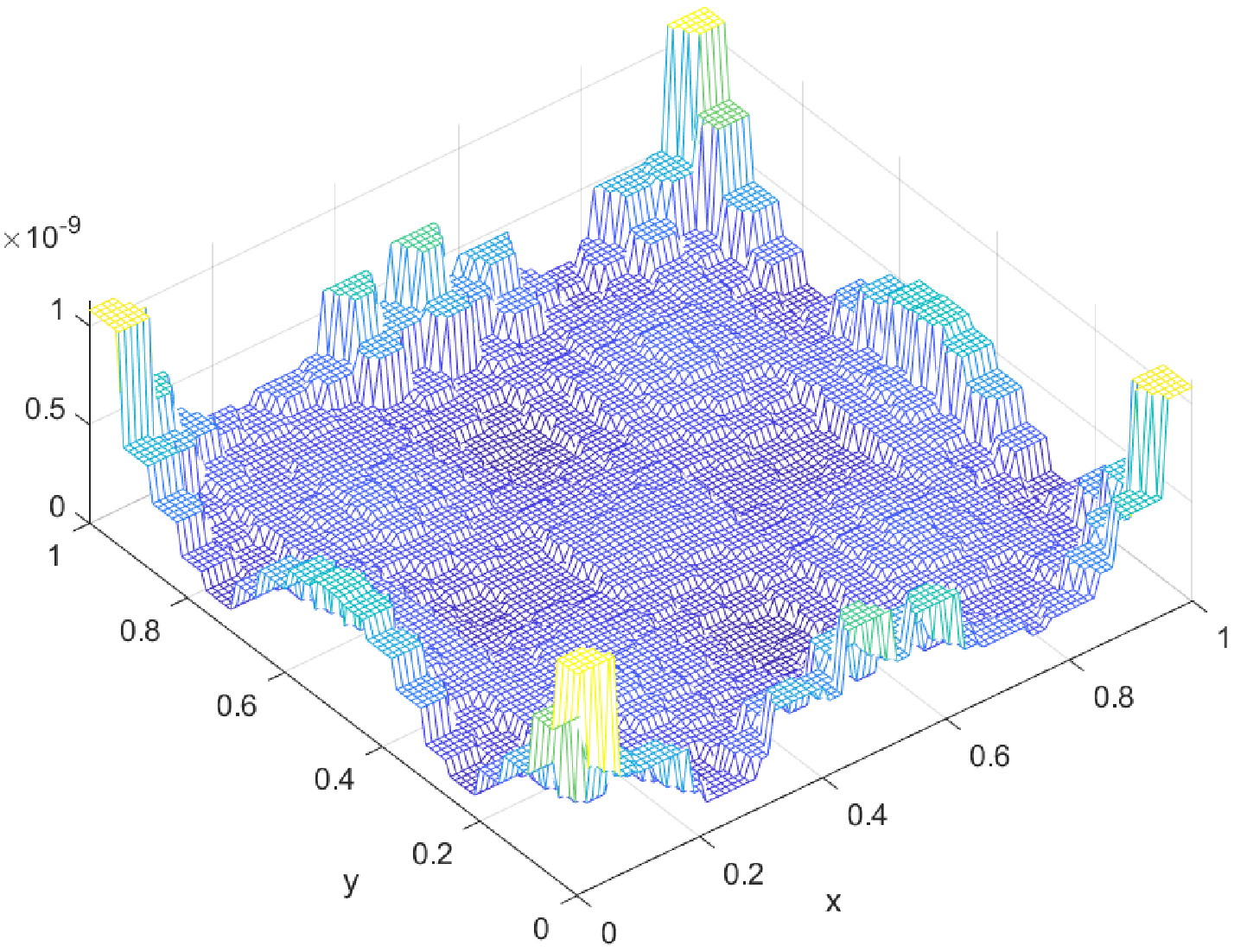}}\vspace*{-4pt}\centerline{\scriptsize Exact solution with four humps}}
}
\vspace*{-8pt}\caption{The $L^2$ error on each element for the two additional tests based on mesh sequence
  $\cT^1_h$ at level $n = 18$ with approximation index $r=5$ }\label{fig:error_boundary}
\end{figure}

We performed two additional tests, with the $L^2$ error on each element shown in
Figure~\ref{fig:error_boundary}.  For the first additional test, we solved the same problem on the
domain $[0,2]\times[-1,1]$ using a mesh given by reflecting the original $n=18$ mesh with respect to
$x=1$, and then reflecting this with respect to $y=0$. This test shows that when the original boundary
elements are moved to the interior of the domain, we still observe the same larger error.  For the
second additional test, we solved the problem on the unit square domain with the original mesh, but
we set the exact solution to be $u(x_1,x_2) = \sin(2\pi x_1)\sin(2\pi x_2)$, which has four humps.
From the figure, we see that the solution is better approximated in the interior where hexagons are
used versus the approximation near the boundary.

\begin{figure}[ht]
\centerline{
\parbox{.5\linewidth}{\centerline{\includegraphics[width=0.60\linewidth]{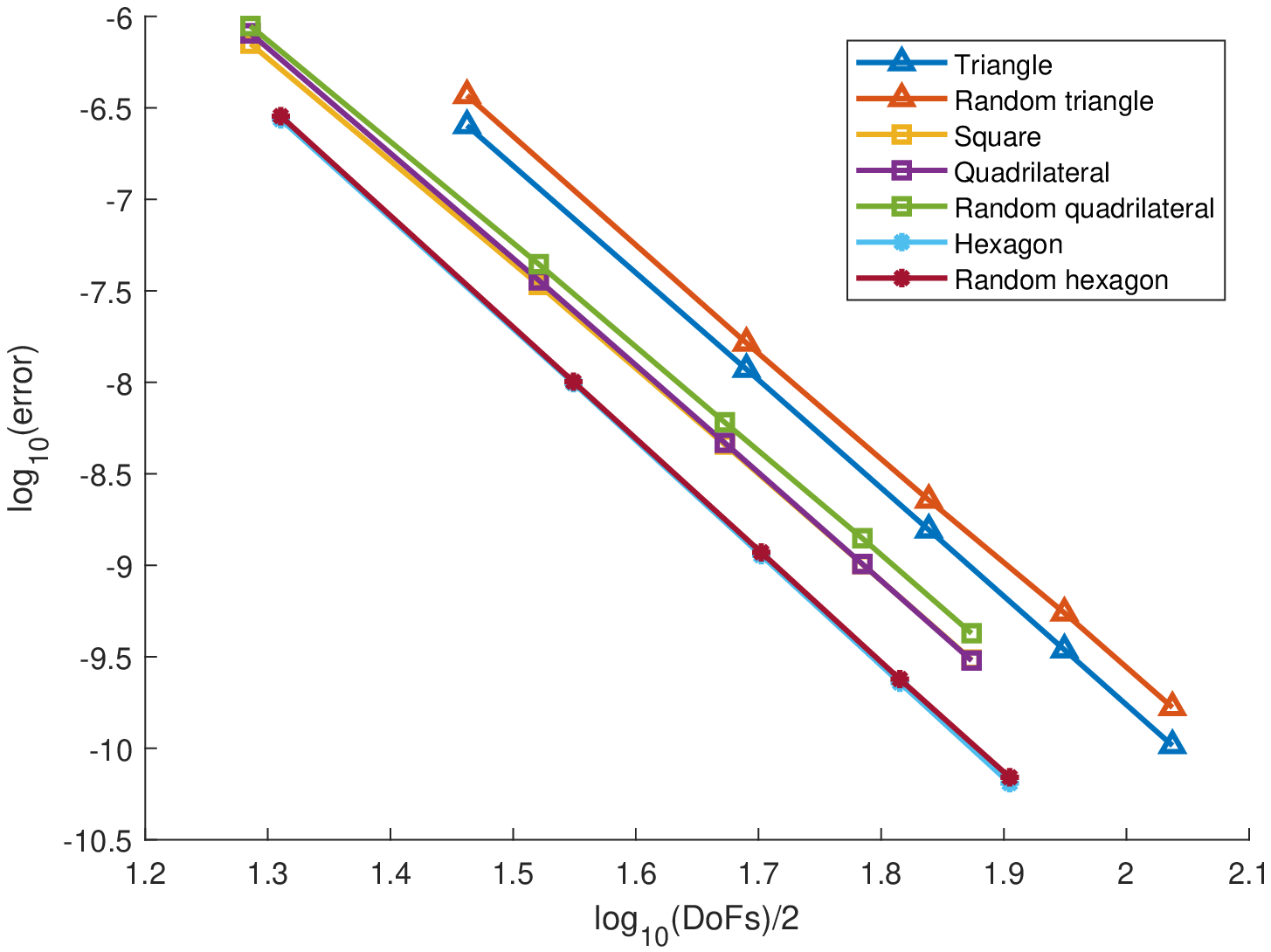}%
\llap{\raisebox{20pt}{\scriptsize $L^2$-norm errors\hspace*{40pt}}}}}
\parbox{.5\linewidth}{\centerline{\includegraphics[width=0.60\linewidth]{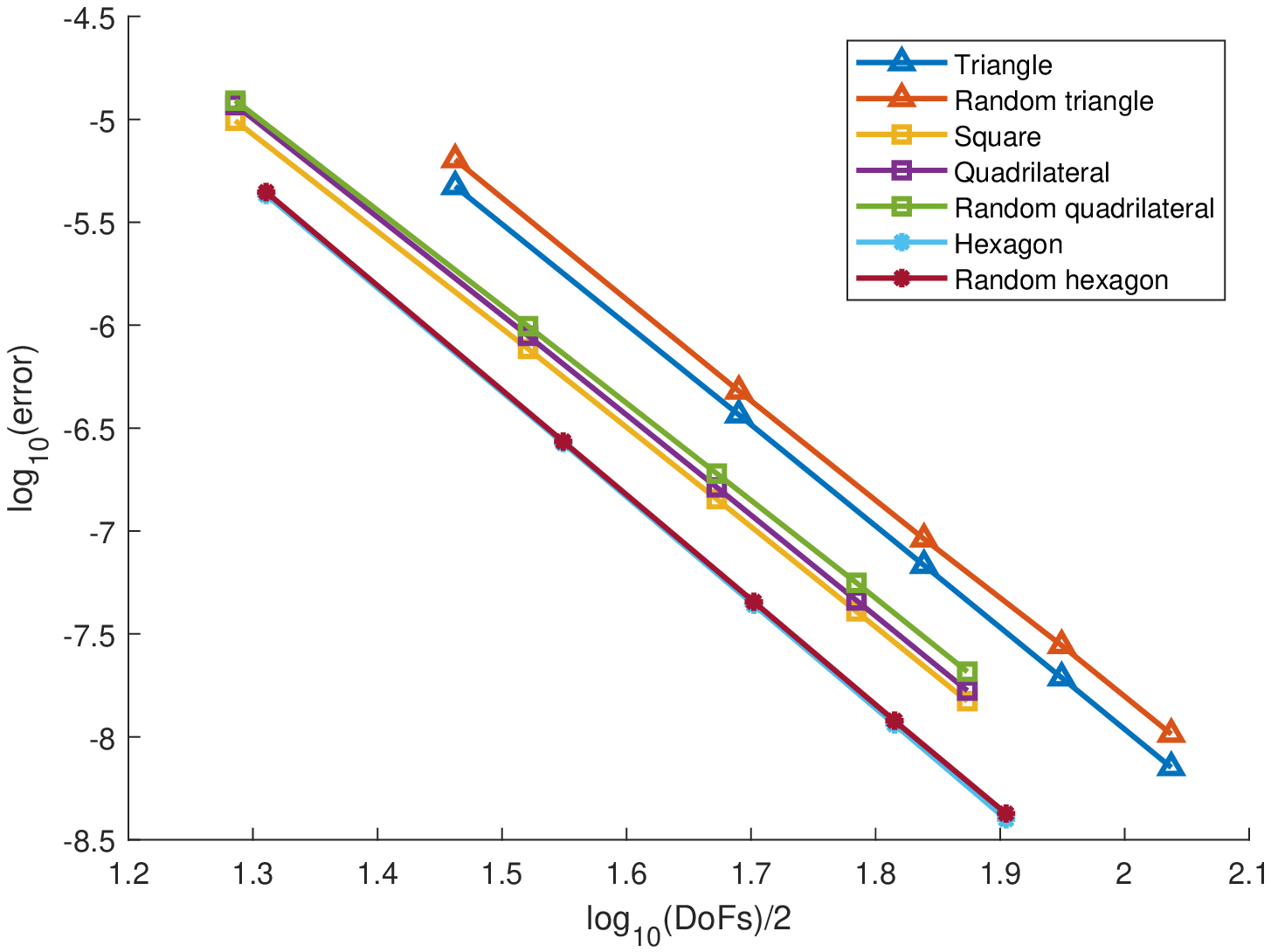}%
\llap{\raisebox{20pt}{\scriptsize $H^1$-seminorm\hspace*{45pt}}}}}
}
\vspace*{-4pt}\caption{Log of the $L^2$-norm and $H^1$-seminorm errors versus half the log of the number of DoFs
  on seven different mesh sequences with $n=6,10,14,18,22$ and $r=5$.}\label{fig:error_dofs}
\end{figure}

To further verify that hexagons are better at approximation, we performed experiments for index
$r=2,3,4,5$ at levels $n=6,10,14,18,22$ on seven different meshes, each emphasizing a fixed number
of edges $N$ for the elements. The first mesh consists of isosceles right triangles, and we distort
it with random noise to get the second mesh. The third mesh consists of squares, the fourth mesh is
a mesh of identical trapezoids, and the fifth mesh consists of quadrilaterals obtained by randomly
distorting the vertices of a square mesh.  The sixth mesh is $\cT_h^1$ (mostly hexagons), and we
distort it with some randomness to get the seventh mesh.  To simplify the presentation, we only show
results for $r=5$ in Figure~\ref{fig:error_dofs}, since the others are similar. We plot the log of
error versus half the log of the number of degrees of freedom for each mesh sequence.  We see that
for the same number of degrees of freedom, hexagonal elements give the best results, followed by
quadrilaterals, with triangular elements giving the worst performance.

\begin{table}[ht]
\caption{Errors and convergence rates for $\cDS_r$ on $\cT^2_h$ meshes.}\label{tab:t3}
\vspace*{-4pt}\begin{center}\small
\begin{tabular}{c|cc|cc|cc|cc}
\hline & \multicolumn{2}{c|}{$r=2$} & \multicolumn{2}{c|}{$r=3$}
& \multicolumn{2}{c|}{$r=4$} & \multicolumn{2}{c}{$r=5$} \\
$n$ & error & rate & error & rate & error & rate & error & rate \\
\hline
\multicolumn{9}{c}{$L^2$-errors and convergence rates}\\
\hline
10 & 2.160e-04 & 3.45 & 8.859e-06 & 4.34 & 3.467e-07 & 5.69 & 1.133e-08 & 6.97 \\
14 & 7.329e-05 & 3.16 & 2.175e-06 & 4.11 & 5.644e-08 & 5.31 & 1.202e-09 & 6.57 \\
18 & 3.452e-05 & 2.95 & 7.927e-07 & 3.96 & 1.530e-08 & 5.12 & 4.376e-10 & 3.97 \\
22 & 1.863e-05 &	3.47 &	3.555e-07	& 4.51 &	5.314e-09 &	5.95 &	8.905e-11 &	8.95 \\
\hline
\multicolumn{9}{c}{$H^1$-seminorm errors and convergence rates}\\
\hline
10 & 3.561e-03 & 2.32 & 1.933e-04 & 3.13 & 8.530e-06 & 4.55 & 3.103e-07 & 5.73 \\
14 & 1.683e-03 & 2.19 & 6.724e-05 & 3.09 & 1.973e-06 & 4.29 & 4.625e-08 & 5.57 \\
18 & 1.018e-03 & 1.97 & 3.144e-05 & 2.98 & 6.952e-07 & 	4.09 & 2.646e-08 & 2.19 \\
22 & 6.712e-04	& 2.34 &	1.730e-05 &	3.36 &	2.969e-07 &	4.78 &	5.973e-09 &	8.37 \\ 

\hline
\end{tabular}
\end{center}
\end{table}


\subsubsection{Not so shape regular meshes of mostly hexagons, $\cT_h^2$}\label{sec:DS_nonregularMesh}

Table~\ref{tab:t3} presents the errors and orders of convergence for the mesh sequence $\cT_h^2$
generated by $n^2$ random initial seeds. We see that the convergence rates are generally correct,
but they are not steady due to the randomness inherent in the mesh refinement process. Of particular
concern are the rates for $n = 18,22$, especially as $r$ increases. We attribute this behavior to the 
poor shape regularity of these two random meshes (recall Table~\ref{tab:chunk}).

\begin{figure}[ht]
\centerline{
\parbox{.5\linewidth}{\centerline{\includegraphics[width=0.60\linewidth]{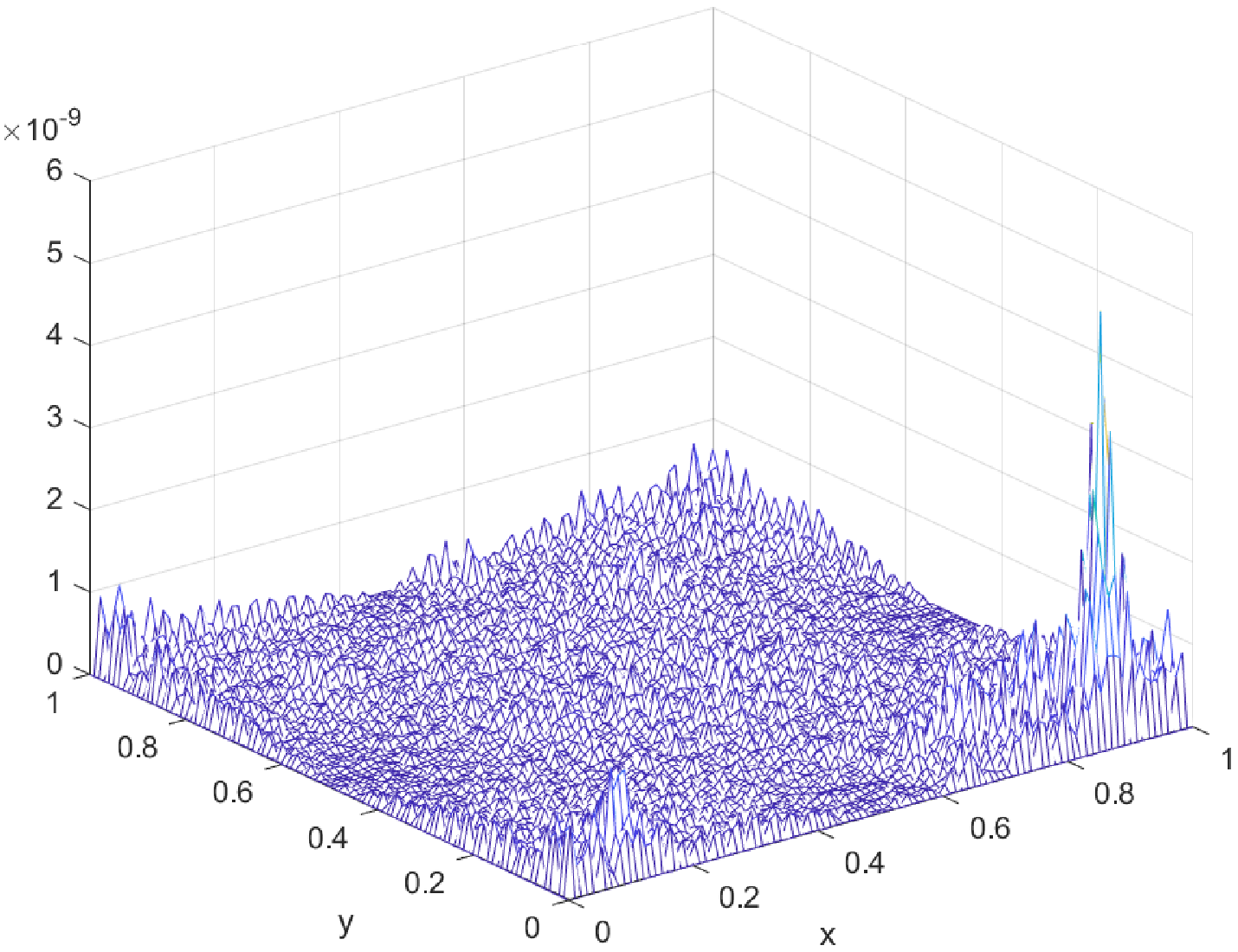}}\centerline{\scriptsize Pointwise error for original $\cT^2_h$}}
\parbox{.5\linewidth}{\centerline{\includegraphics[width=0.60\linewidth]{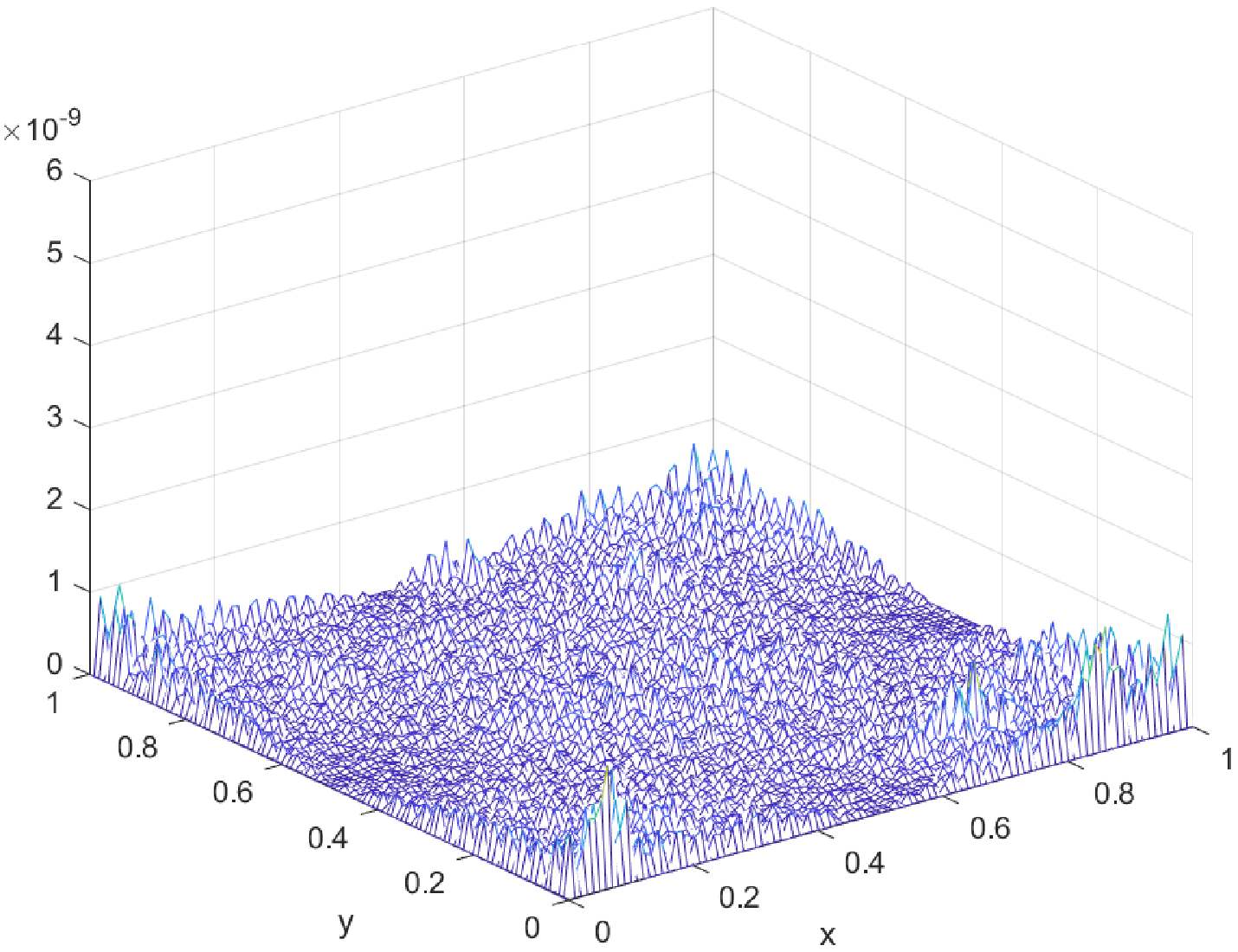}}\centerline{\scriptsize Pointwise error for modified $\cT^2_h$}}
}
\vspace*{-4pt}\caption{The plot of pointwise error 
  for $\cT^2_h$ with level $n = 18$ and approximation index $r=5$ before and after modifying the
  mesh. The two small edges were removed by removing vertices located at $(0.108,0.050)$ and
  $(0.890,0.057)$.}\label{fig:error_t3}
\end{figure}

An examination of the spatial distribution of the error for $n = 18$, as shown on the left in
Figure~\ref{fig:error_t3}, suggests that the error is exceptionally large near one corner.  The
$n = 18$ mesh has two edges that are relatively very short containing the vertices $(0.108,0.050)$
and $(0.890,0.057)$, and the $n = 22$ mesh has five short edges.  We created the modified $\cT^2_h$
meshes by removing one vertex of each short edge.  As can be seen in Table~\ref{tab:chunk}, the
shape regularity parameters of the elements of the modified mesh are more uniform.  The right plot
in Figure~\ref{fig:error_t3} shows that the error is reduced without the offending edges. The overall
error and convergence results for the modified mesh are presented in Table~\ref{tab:t3-2}, and they
are closer to the expected rates.

\begin{table}[ht]
\caption{Errors and convergence rates for $\cDS_r$ on modified $\cT^2_h$ meshes.}\label{tab:t3-2}
\vspace*{-4pt}\begin{center}\small
\begin{tabular}{c|cc|cc|cc|cc}
\hline & \multicolumn{2}{c|}{$r=2$} & \multicolumn{2}{c|}{$r=3$}
& \multicolumn{2}{c|}{$r=4$} & \multicolumn{2}{c}{$r=5$} \\
\hline
$n$ & error & rate & error & rate & error & rate & error & rate \\
\hline
\multicolumn{9}{c}{$L^2$-errors and convergence rates}\\
\hline
18 & 3.454e-05 & 3.30 & 8.172e-07 & 4.29 & 1.544e-08 & 5.68 & 3.080e-10 & 5.97 \\
22 & 1.8811e-05	& 3.26 &	3.6047e-07 & 4.39 & 5.4763e-09 & 5.56 &	8.1513e-11 & 7.13 
  \\
\hline
\multicolumn{9}{c}{$H^1$-seminorm errors and convergence rates}\\
\hline
18 & 1.018e-03 & 2.20 & 3.194e-05 & 3.26 & 6.992e-07 & 4.55 & 1.553e-08 & 4.78 \\
22 & 6.7623e-04	& 2.19 & 1.7434e-05 & 3.25 & 3.0345e-07 & 4.48 &	4.9949e-09 & 6.09 \\
\hline
\end{tabular}
\end{center}
\end{table}


\subsection{Direct mixed spaces}

We now consider the direct mixed finite elements $\V_r^s\times W_s$ derived in
Section~\ref{sec:implementMixed}. These are implemented both in hybrid form
(Section~\ref{sec:hybrid}) and as $H(\Div)$-conforming elements (Section~\ref{sec:conforming}),
which, of course, provide the same results.

The $L^2$ and $H^1$-seminorm errors and convergence orders for the mesh sequence $\cT_h^1$ with
$r=(0,\!)\,1,2,3$ appear in Tables~\ref{tab:at1-h2-r}--\ref{tab:at1-h2-f}.  The theory predicts that
the scalar $p$, the vector $\u$, and the divergence $\div\u$ should attain the order of
approximation $s+1$, $r+1$, and $s+1$, respectively, for the reduced ($s=r-1$) and full ($s=r$)
$H(\Div)$-approximation spaces. We see rates of convergence that are close to the theoretical ones.
Moreover, the errors for $\cT_h^1$ are a bit smaller than what we see for meshes of trapezoids, due
to having many elements with more than four edges.

\begin{table}[ht]
\caption{Errors and convergence rates in $L^2$ for
direct mixed reduced $H(\Div)$-approximation on $\cT_h^1$ meshes.}\label{tab:at1-h2-r}
\vspace*{-4pt}\begin{center}\small
\begin{tabular}{c|cc|cc|cc}
\hline
&  \multicolumn{2}{c|}{$\|p-p_h\|$} & \multicolumn{2}{c|}{$\|\u-\u_h\|$}
& \multicolumn{2}{c}{$\|\div(\u-\u_h)\|$}\\
$n$ & error & rate & error & rate & error & rate\\
\hline
\multicolumn{7}{c}{$r=1$, reduced $H(\Div)$-approximation}\\
\hline
10 & 1.308e-01 & 1.10 & 1.820e-02 & 2.05 & 1.277e-01 & 1.02 \\
14 & 9.196e-02 & 1.05 & 9.199e-03 & 2.03 & 9.084e-02 & 1.01 \\
18 & 7.104e-02 & 1.03 & 5.539e-03 & 2.02 & 7.051e-02 & 1.01 \\
22 & 5.791e-02 & 1.02 & 3.698e-03 & 2.01 & 5.763e-02 & 1.01 \\

\hline
\multicolumn{7}{c}{$r=2$, reduced $H(\Div)$-approximation}\\
\hline
10 & 8.640e-03 & 2.04 & 5.053e-04 & 3.04 & 8.639e-03 & 2.04 \\
14 & 4.363e-03 & 2.03 & 1.825e-04 & 3.03 & 4.363e-03 & 2.03 \\
18 & 2.624e-03 & 2.02 & 8.545e-05 & 3.02 & 2.624e-03 & 2.02 \\
22 & 1.750e-03 & 2.02 & 4.666e-05 & 3.01 & 1.750e-03 & 2.02 \\

\hline
\multicolumn{7}{c}{$r=3$, reduced $H(\Div)$-approximation}\\
\hline
10 & 3.858e-04 & 3.07 & 1.831e-05 & 4.06 & 3.858e-04 & 3.07 \\
14 & 1.385e-04 & 3.05 & 4.710e-06 & 4.04 & 1.385e-04 & 3.05 \\
18 & 6.464e-05 & 3.03 & 1.713e-06 & 4.02 & 6.464e-05 & 3.03 \\
22 & 3.522e-05 & 3.03 & 7.643e-07 & 4.02 & 3.522e-05 & 3.03 \\

\hline
\end{tabular}
\end{center}
\end{table}

\begin{table}[ht]
\caption{Errors and convergence rates in $L^2$ for
direct mixed full $H(\Div)$-approximation on $\cT_h^1$ meshes.}\label{tab:at1-h2-f}
\vspace*{-4pt}\begin{center}\small
\begin{tabular}{c|cc|cc|cc}
\hline
&  \multicolumn{2}{c|}{$\|p-p_h\|$} & \multicolumn{2}{c|}{$\|\u-\u_h\|$}
& \multicolumn{2}{c}{$\|\div(\u-\u_h)\|$}\\
$n$ & error & rate & error & rate & error & rate\\
\hline
\multicolumn{7}{c}{$r=0$, full $H(\Div)$-approximation}\\
\hline

10 & 1.299e-01 & 1.07 & 6.167e-02 & 1.36 & 1.277e-01 & 1.02 \\
14 & 9.170e-02 & 1.04 & 3.970e-02 & 1.31 & 9.084e-02 & 1.01 \\
18 & 7.093e-02 & 1.02 & 2.883e-02 & 1.27 & 7.051e-02 & 1.01 \\
22 & 5.786e-02 & 1.01 & 2.245e-02 & 1.25 & 5.763e-02 & 1.01 \\

\hline
\multicolumn{7}{c}{$r=1$, full $H(\Div)$-approximation}\\
\hline
10 & 8.641e-03 & 2.04 & 2.403e-03 & 2.38 & 8.639e-03 & 2.04 \\
14 & 4.363e-03 & 2.03 & 1.094e-03 & 2.34 & 4.363e-03 & 2.03 \\
18 & 2.624e-03 & 2.02 & 6.133e-04 & 2.30 & 2.624e-03 & 2.02 \\
22 & 1.759e-03 & 1.99 & 3.888e-04 & 2.27 & 1.750e-03 & 2.02 \\

\hline
\multicolumn{7}{c}{$r=2$, full $H(\Div)$-approximation}\\
\hline
10 & 3.858e-04 & 3.07 & 7.535e-05 & 3.37 & 3.858e-04 & 3.07 \\
14 & 1.385e-04 & 3.05 & 2.420e-05 & 3.38 & 1.385e-04 & 3.05 \\
18 & 6.464e-05 & 3.03 & 1.038e-05 & 3.37 & 6.464e-05 & 3.03 \\
22 & 3.522e-05 & 3.03 & 5.288e-06 & 3.36 & 3.522e-05 & 3.03 \\ 

\hline
\multicolumn{7}{c}{$r=3$, full $H(\Div)$-approximation}\\
\hline

10 & 1.372e-05 & 4.13 & 2.572e-06 & 4.52 & 1.372e-05 & 4.13 \\
14 & 3.459e-06 & 4.10 & 5.879e-07 & 4.39 & 3.459e-06 & 4.10 \\
18 & 1.243e-06 & 4.07 & 1.987e-07 & 4.32 & 1.243e-06 & 4.07 \\
22 & 5.502e-07 & 4.06 & 8.451e-08 & 4.26 & 5.502e-07 & 4.06 \\

\hline
\end{tabular}
\end{center}
\end{table}

The errors and orders of convergence of the modified $\cT_h^2$ mesh sequence are given in
Tables~\ref{tab:at1-mh3-r}--\ref{tab:at1-mh3-f}.  We see the expected results.

\begin{table}[ht]
\caption{Errors and convergence rates in $L^2$ for
		direct mixed reduced $H(\Div)$-approx.\ on modified $\cT_h^2$ meshes.}\label{tab:at1-mh3-r}
	\vspace*{-12pt}\begin{center}\small
		\begin{tabular}{c|cc|cc|cc}
			\hline
			&  \multicolumn{2}{c|}{$\|p-p_h\|$} & \multicolumn{2}{c|}{$\|\u-\u_h\|$}
			& \multicolumn{2}{c}{$\|\div(\u-\u_h)\|$}\\
			$n$ & error & rate & error & rate & error & rate\\
			\hline
			\multicolumn{7}{c}{$r=1$, reduced $H(\Div)$-approximation}\\
			\hline
10 & 1.290e-01 & 1.24 & 1.770e-02 & 2.29 & 1.260e-01 & 1.15 \\
14 & 9.109e-02 & 1.02 & 8.997e-03 & 1.98 & 9.001e-02 & 0.98 \\
18 & 7.039e-02 & 1.13 & 5.429e-03 & 2.21 & 6.988e-02 & 1.11 \\
22 & 5.734e-02 & 1.10 & 3.619e-03 & 2.18 & 5.707e-02 & 1.09 \\
			\hline
			\multicolumn{7}{c}{$r=2$, reduced $H(\Div)$-approximation}\\
			\hline
10 & 8.635e-03 & 2.23 & 5.013e-04 & 3.24 & 8.634e-03 & 2.23 \\
14 & 4.308e-03 & 2.04 & 1.785e-04 & 3.02 & 4.308e-03 & 2.03 \\
18 & 2.616e-03 & 2.19 & 8.487e-05 & 3.26 & 2.616e-03 & 2.19 \\
22 & 1.719e-03 & 2.25 & 4.649e-05 & 3.23 & 1.719e-03 & 2.25 \\
			\hline
			\multicolumn{7}{c}{$r=3$, reduced $H(\Div)$-approximation}\\
			\hline
10 & 3.878e-04 & 3.38 & 1.992e-05 & 4.37 & 3.878e-04 & 3.38 \\
14 & 1.384e-04 & 3.02 & 5.102e-06 & 3.99 & 1.384e-04 & 3.02 \\
18 & 6.516e-05 & 3.30 & 1.889e-06 & 4.36 & 6.516e-05 & 3.30 \\
22 & 3.514e-05 & 3.31 & 8.363e-07 & 4.37 & 3.514e-05 & 3.31 \\
			\hline
		\end{tabular}
	\end{center}
\end{table}

\begin{table}[ht]
	\caption{Errors and convergence rates in $L^2$ for
		direct mixed full $H(\Div)$-approx.\ on modified $\cT_h^2$ meshes.}\label{tab:at1-mh3-f}
	\vspace*{-4pt}\begin{center}\small
		\begin{tabular}{c|cc|cc|cc}
			\hline
			&  \multicolumn{2}{c|}{$\|p-p_h\|$} & \multicolumn{2}{c|}{$\|\u-\u_h\|$}
			& \multicolumn{2}{c}{$\|\div(\u-\u_h)\|$}\\
			$n$ & error & rate & error & rate & error & rate\\
		\hline
	\multicolumn{7}{c}{$r=0$, full $H(\Div)$-approximation}\\
	\hline		
10 & 1.282e-01 & 1.20 & 5.915e-02 & 1.59 & 1.260e-01 & 1.15 \\
14 & 9.089e-02 & 1.01 & 3.577e-02 & 1.47 & 9.001e-02 & 0.98 \\
18 & 7.030e-02 & 1.13 & 2.701e-02 & 1.23 & 6.988e-02 & 1.11 \\
22 & 5.730e-02 & 1.10 & 2.005e-02 & 1.60 & 5.707e-02 & 1.09  \\
			\hline
			\multicolumn{7}{c}{$r=1$, full $H(\Div)$-approximation}\\
			\hline
10 & 8.635e-03 & 2.23 & 1.892e-03 & 2.67 & 8.634e-03 & 2.23 \\
14 & 4.308e-03 & 2.04 & 8.562e-04 & 2.32 & 4.308e-03 & 2.03 \\
18 & 2.616e-03 & 2.19 & 4.903e-04 & 2.44 & 2.616e-03 & 2.19 \\
22 & 1.719e-03 & 2.25 & 3.142e-04 & 2.39 & 1.719e-03 & 2.25  \\
			\hline
			\multicolumn{7}{c}{$r=2$, full $H(\Div)$-approximation}\\
			\hline
10 & 3.881e-04 & 3.38 & 6.546e-05 & 3.69 & 3.881e-04 & 3.38 \\
14 & 1.384e-04 & 3.02 & 1.945e-05 & 3.55 & 1.384e-04 & 3.02 \\
18 & 6.516e-05 & 3.30 & 8.982e-06 & 3.39 & 6.516e-05 & 3.30 \\
22 & 3.514e-05 & 3.31 & 4.448e-06 & 3.77 & 3.514e-05 & 3.31  \\
			\hline
			\multicolumn{7}{c}{$r=3$, full $H(\Div)$-approximation}\\
			\hline
10 & 1.299e-05 & 4.59 & 2.473e-06 & 5.15 & 1.299e-05 & 4.59 \\
14 & 3.270e-06 & 4.04 & 5.434e-07 & 4.44 & 3.270e-06 & 4.04 \\
18 & 1.188e-06 & 4.44 & 2.220e-07 & 3.92 & 1.188e-06 & 4.44 \\
22 & 5.259e-07 & 4.37 & 1.021e-07 & 4.17 & 5.259e-07 & 4.37 \\
			\hline
		\end{tabular}
	\end{center}
\end{table}


\section{Summary and Conclusions}\label{sec:conc}

We defined direct serendipity finite elements on general closed, nondegenerate, and convex polygons
$E_N$ with $N$ vertices for any index of approximation~$r$. A direct serendipity element has its
function space of the form of polynomials plus supplemental functions, i.e.,
\begin{equation}\label{eq:serendipityFormConc}
  \cDS_r(E_N) = \Po_r(E_N)\oplus\Supp_r^\cDS(E_N),\quad r\geq1,
\end{equation}
with the supplemental space $\Supp_r^\cDS(E_N)$ being of minimal local dimension subject to the
requirement of global $H^1$-conformity. For higher order finite element spaces with $r\ge N-2$, the
supplemental space $\Supp_r^\cDS(E_N)$ has dimension $\frac{1}{2}N(N-3)$, which is the number of
pairs of nonadjacent edges. This fact inspires our construction \eqref{eq:supplementSpace}, for
which different choices of $\lambda_{i,j}$ and $R_{i,j}$ give rise to different spaces. Each index
${i,j}$ represents a pair of nonadjacent edges $e_i$ and $e_j$ of $E_N$.
Simple choices for $\lambda_{i,j}$ and $R_{i,j}$ can be made, as given in \eqref{eq:simpleLambdaHV}
and \eqref{eq:R-simple}. The lower order direct serendipity finite element spaces with $r<N-2$, are
given as the subset of functions in $\cDS_{N-2}(E_N)$ that restrict to polynomials of degree $r$
on $\d E_N$.  Taking nodal DoFs, we constructed nodal bases for the direct serendipity spaces.

By the de Rham theory, each direct serendipity element $\cDS_{r+1}(E_N)$ gives rise to a reduced and
a full direct mixed $H(\Div)$-approximation mixed finite element
\begin{align}
\label{eq:generalMixedRedSumm}
\V_r^{r-1}(E_N) &= \Curl\,\cDS_{r+1}(E_N)\oplus\x\Po_{r-1}(E_N)
= \Po_r^2(E_N)\oplus\Supp_r^\V(E_N),\quad r\ge1,\\
\label{eq:generalMixedSumm}
\V_r^{r}(E) &= \Curl\,\cDS_{r+1}(E_N)\oplus\x\Po_{r}(E_N)\\\nonumber
&\qquad\qquad\qquad\qquad\quad\ 
                  = \Po_r^2(E_N)\oplus\x\tilde\Po_{r}(E_N)\oplus\Supp_r^\V(E_N),\quad r\ge0,
\end{align}
respectively, where $\Supp_r^\V(E) = \Curl\,\Supp_{r+1}^\cDS(E_N)$ has minimal local dimension
subject to the requirement of global $H(\Div)$-conformity.  These mixed elements can be implemented
globally in the hybrid form of the mixed method without the need of a global basis.  However, we
also provided an explicit conforming global basis that we constructed locally on each $E_N$ using
the basis of $\cDS_{r+1}(E_N)$.

The convergence theory handled the polygonal geometry through a continuous dependence argument over
a compact set of perturbations.  Assuming that the meshes are shape regular as $h\to0$
(Definition~\ref{defn:shape_regular}) and that the functions $\lambda_{i,j}$ and $R_{i,j}$ in
\eqref{eq:supplementSpace} are chosen to be continuously differentiable with respect to the vertices
of the element (i.e., Assumption~\ref{assumption:niceDS}), we obtained optimal approximation rates
for the elements in Theorems~\ref{thm:approxDS} and~\ref{thm:brambleMixed}.

We presented and discussed numerical results from finite element numerical solutions of Poisson's
equation. The convergence rates were consistent with the theory, Theorem~\ref{thm:convergence}, and
provided confirmation of the optimal order of accuracy of the finite element approximations.  We
found that mesh shape regularity was quite important in terms of the observed error.  In particular,
we found that short edges, which lead to a poor (i.e., small) shape regularity parameter, could also
result in a poor approximation in that region of the mesh.  Removing such edges greatly improved the
approximation and convergence rates.  We also observed that meshes that emphasize elements with many
edges per element out perform meshes with fewer edges per element. This observation, as well as the
need for flexible meshing in some applications, can be considered justification for using polygonal
elements.



\end{document}